\documentclass[11pt]{amsart}
\usepackage{latexsym, amsmath, amssymb}
\usepackage{epsfig}
\usepackage{amscd}
\usepackage{pstricks,pst-node,cite}

\ifx\blackandwhite\undefined
  \psset{linecolor=black} \else
  \psset{linecolor=black}
\fi

\newcommand{\middlearrow}{\lput{:U}{\pspicture[shift=0](0,0)(0,0)
\psline[arrows=->,arrowscale=1.5](2.2pt,0)(2.3pt,0)\endpspicture}}

\newtheorem{thm}{Theorem}[section]
\newtheorem{lem}[thm]{Lemma}

\newtheorem{conj}[thm]{Conjecture}
\newtheorem{rem}[thm]{Remark}

\begin{document}

\title[The quantum $\mathfrak{sl}(n,\mathbb{C})$ representation theory and its applications]
{The quantum $\mathfrak{sl}(n,\mathbb{C})$ representation theory and its applications}

\author{Myeong--Ju Jeong}
\address{Korea Science Academy, Busanjin Gu, Busan, 614-822, Korea}
\email{mjjeong@kaist.ac.kr}

\author{Dongseok Kim}
\address{Department of Mathematics \\Kyonggi University
\\ Suwon, 443-760, Korea}
\email{dongseok@kgu.ac.kr}
\thanks{Corresponding author. Tel:+82-31-249-9613; FAX
number:+82-31-253-1165. This work was supported by the Korea Science and Engineering Foundation
(KOSEF) grant funded by the Korea government(MEST) (2009-0071547).}

\begin{abstract}
In this paper, we study the quantum $\mathfrak{sl}(n)$ representation category using the web space.
Specially, we extend $\mathfrak{sl}(n)$ web space for $n\ge 4$ as generalized Temperley-Lieb algebras.
As an application of our study, we find that the HOMFLY polynomial $P_n(q)$ specialized to a one variable polynomial can be computed
by a linear expansion with respect to a presentation of the quantum representation category of $\mathfrak{sl}(n)$.
Moreover, we correct the false conjecture~\cite{PS:superiod} given by Chbili,
which addresses the relation between some link polynomials of a periodic link
and its factor link such as Alexander polynomial $(n = 0)$ and Jones polynomial $(n = 2)$ and prove
the corrected conjecture not only for HOMFLY polynomial but also for the colored HOMFLY polynomial specialized to a one variable polynomial.
\end{abstract}


\maketitle

\section{Introduction}

The discovery of the Jones polynomial~\cite{Jones:subfactor, Jones:braid} brought
a Renaissance of knot theory and its generalizations have been studied in many
different ways \cite{CK:homology, FK:canonical, Khovanov:colored, Kuperberg:spiders, kirbyMelvin:witten,
Murakami:coloredjones, Vybornov:yang, Witten:pathint}. Using the representation theory of
complex simple Lie algebras, Reshetikhin and Turaev found quantized simple
Lie algebras invariants of links and 3-manifolds~\cite{RT:1,RT:2} and these invariants
have been studied extensively~\cite{chbili, chbili:qm, Khovanov:sl3, OY:quantum,
KR:factor, MOY:Homfly, yokota:skeinforn}.

In the present paper, we study the quantum $\mathfrak{sl}(n)$ representation theory related to
the HOMFLY polynomials of periodic links.
Murasugi found a strong relation between the Alexander polynomials
of a periodic link $L$ and its factor link $\overline L$~\cite{mu:alexander} and a similar relation
for the Jones polynomials~\cite{mu:jones}. There are various results to decide the periodicity of links
\cite{Przytycki:criterion, traczyk:period3, JP5, yokota:skein, yokota:Jones, yokota:kauffman}.
A conjecture for the relation between HOMFLY polynomials $P_n(q)$ specialized to a one variable polynomial of a
periodic link $L$ and its factor link $\overline L$ was found as follows~\cite{chbili}.

\begin{conj} [\cite{chbili}]
Let $p$ be a positive integer and $L$ be a $p$-periodic link in $\mathbb{S}^3$ with its factor link $\overline{L}$.
Then,
$$P_n(L) \equiv P_n(\overline{L})^p \hskip .7cm modulo \hskip .2cm
\mathcal{A}_n,$$ where $\mathcal{A}_n$ is the ideal of $\mathbb{Z}[q^{\pm \frac 12}]$
generated by $p$ and $[n]^p-[n]$. \label{conj1}
\end{conj}
\noindent The quantum integers are defined as
$$
[n]=\frac{q^{\frac{n}{2}}-q^{-\frac{n}{2}}}{q^{\frac{1}{2}}-q^{-\frac{1}{2}}},\quad
 [n]!= [n][n-1]\ldots [2][1], \quad
\left[\begin{matrix} n \\ k \end{matrix}\right] =
\frac{[n]!}{[n-k]![k]!}.
$$

For Conjecture~\ref{conj1}, Chbili provided a proof for $n=3$ using the representation theory of the quantum
$\mathfrak{sl}(3)$ \cite{chbili}. There were subsequent studies on the conjecture~\cite{CL:period}.
But, it was shown that \textit{Conjecture~\ref{conj1} is \textbf{false} for
$n\ge 4$}~\cite{PS:superiod}. The counterexamples show that even if the link is
colored by the vector representation, the congruence can be involved with other fundamental
representations. Focused on the quantum representation category of
$\mathfrak{sl}(n,\mathbb{C})$, we may modify the original conjecture as follows.

\begin{conj}
Let $p$ be a positive integer and $L$ be a $p$-periodic link in $\mathbb{S}^3$
with its factor link $\overline{L}$. Then,
$$P_n(L) \equiv P_n(\overline{L})^p  \hskip .7cm modulo \hskip .2cm
\mathcal{I}_n,$$ where $\mathcal{I}_n$ is the ideal of
$\mathbb{Z}[q^{\pm \frac 12}]$ generated by $p$ and
$\left[\begin{matrix}n\\
i\end{matrix}\right]^p-\left[\begin{matrix}n\\
i\end{matrix}\right]$ for $i=1, 2, \ldots, \lfloor \frac{n}{2}
\rfloor$. \label{modifiedconj}
\end{conj}

The study of a presentation of the quantum representation category of $\mathfrak{sl}(n)$ leads
us to a powerful computation method of the HOMFLY polynomial $P_n(q)$ specialized to a one variable polynomial, a linear expansion of webs, and its generalization to
the colored $\mathfrak{sl}(n)$ HOMFLY polynomial $G_n(L,\mu)$ specialized to a one variable polynomial.
Then we not only prove Conjecture~\ref{modifiedconj} in Theorem~\ref{mainthm} but also we show the following theorem that
Conjecture~\ref{modifiedconj} remains true for $G_n(L,\mu)$.

\begin{thm}
Let $p$ be a positive integer and $L$ be a $p-$periodic link in
$\mathbb{S}^3$ with its factor link $\overline{L}$. Let $\mu$ be a
$p$-periodic coloring of $L$ and $\overline{\mu}$ be the induced
coloring of $\overline{L}$. Then for $n\ge 0$,

$$G_n(L,\mu) \equiv G_n(\overline{L},\overline{\mu})^p \hskip .7cm modulo \hskip .2cm
\mathcal{I}_n,$$ where $\mathcal{I}_n$ is the ideal of
$\mathbb{Z}[q^{\pm \frac 12}]$ generated by $p$ and
$\left[\begin{matrix}n\\
i\end{matrix}\right]^p-\left[\begin{matrix}n\\
i\end{matrix}\right]$ for $i=1, 2, \ldots, \lfloor \frac{n}{2}
\rfloor$. \label{mainthm4}
\end{thm}

Recently, there are significant progresses on the representation theory of
the quantum $\mathfrak{sl}(n)$~\cite{Morrison, SW:graph, math.QA/0601209}. In particular, a complete set of relations
for the representation theory of
the quantum $\mathfrak{sl}(n)$ which contain
our relations in Figure~\ref{firstrel}, Lemma~\ref{rectangular} and new relation called
'Kekul$\acute{e}$ relation', which was first found by the second author
for $\mathfrak{sl}(4)$~\cite{Dongseok:thesis},
is conjectured in~\cite{Morrison}.  Furthermore, it was also proven that Remark~\ref{remcon} is false~\cite{Morrison}.

The outline of this paper is as follows. In section~\ref{Review},
we review the HOMFLY polynomials and the colored HOMFLY polynomials specialized to a one variable polynomial.
In section~\ref{relation}, we develop the representation theory of
the quantum $\mathfrak{sl}(n)$. We show that the quantum $\mathfrak{sl}(n)$ skein module of the plane or sphere
has dimension $1$ using the relations we have found. In section~\ref{invariants} we prove
Conjecture~\ref{modifiedconj} and show the conjecture holds for the the colored $\mathfrak{sl}(n)$
HOMFLY polynomial specialized to a one variable polynomial. In section~\ref{discussion}, we compare our result with previous works.

\section{The HOMFLY polynomials and the colored HOMFLY polynomials specialized to a one variable polynomial}\label{Review}

A \emph{link} $L$ is a disjoint union of circles embedded in three dimensional sphere $\mathbb{S}^3$,
and a \emph{knot} $K$ is a link with only one component.
Here, we assume all links are PL. A link $L$ in $\mathbb{S}^3$ is $p$-\emph{periodic} if there exists a periodic homeomorphism $h$
of order $p$ such that $fix(h)\cong S^1, h(L)=L$ and $fix(h)\cap L=\emptyset$
where $fix(h)$ is the set of fixed points of $h$. It is well known that if we consider
$\mathbb{S}^3$ as $\mathbb{R}^3 \cup \{\infty\}$, we can assume that $h$ is a rotation by
$2\pi/p$ angle around the $z-$axis. Let $G= \mathbb{Z}/p\mathbb{Z}$ denote
the group of homeomorphisms of $\mathbb{S}^3$ generated by $h$, and let
$\pi$ denote the covering map $\mathbb{S}^3 \rightarrow \mathbb{S}^3 /G$, branched along $z$-axis.
We call $\overline{L} = \pi(L)$ the factor link of $L$.
For other terms and definitions of knot theory, we refer to \cite{Adams}.

Now we define the HOMFLY polynomial specialized to a one variable polynomial.
For the rest of paper, all HOMFLY polynomial and colored HOMFLY polynomial are specialized to a one variable polynomial unless
we state differently.
For a nonnegative integer $n$, the HOMFLY polynomial $P_n(q)$ specialized to a one variable polynomial can be calculated
uniquely by the following skein relations:

\begin{gather}\nonumber P_{n}(\emptyset) =1, \\\nonumber P_{n}(
\begin{pspicture}[shift=-.07](-.17,-.17)(.17,.17) \pscircle(0,0){.15}
\end{pspicture} \cup D)= (\frac{q^{\frac n2}- q^{-\frac n2}}{q^{\frac 12}- q^{-\frac 12}}) P_{n}(D),
 \\\nonumber q^{\frac n2}P_{n}(L_+) - q^{-\frac n2}P_{n}(L_-) = (q^{\frac 12}- q^{-\frac 12}) P_n(L_0), \end{gather}
where $\emptyset$ is the empty diagram,
$\begin{pspicture}[shift=-.07](-.17,-.17)(.17,.17) \pscircle(0,0){.15} \end{pspicture}$
is the trivial knot and $L_+, L_-$ and $L_0$ are skein triple, three diagrams which are identical except at
one crossing as in Figure~\ref{local}.

The HOMFLY polynomial of links can be recovered from the representation
theory of the quantum $\mathfrak{sl}(n)$. For $n=0$, we use the special linear Lie superalgebra $\mathfrak{gl}(1|1)$ to
find that $P_0(q)$ is the Alexander polynomial~\cite{khvanov:comm}. For $n=1$ and for
any link, $P_1(q)=1$. For $n=2$, $P_2(q)$ is the Jones polynomial~\cite{Jones:braid, RT:1, RT:2, Witten:pathint}.
The polynomial $P_n(q)$ can be computed by linearly expanding each crossing
into a sum of diagrams of planar trivalent graphs where the edges of these
planar graphs are oriented and colored by $1$ or $2$ as in Figure~\ref{planar}~\cite{MOY:Homfly}.

To define the colored HOMFLY polynomial, we review the quantum representation category of
$\mathfrak{sl}(n,\mathbb{C})$~\cite{Kuperberg:spiders, KR:factor}.
The color $1$ of an edge in the definition of
HOMFLY polynomial presents the vector representation $V$
of the quantum $\mathfrak{sl}(n)$, and the color $2$ for its exterior power $\wedge^2 V$.
The trivalent vertex is the unique (up to scaling)
intertwiner of $V^{\otimes 2}\otimes \wedge^2 V$.
This setup works for arbitrary exterior powers of $V$~\cite{MOY:Homfly}.
Oriented edges of graphs in their calculus carry colors from $1$ to $n-1$ that denote
the fundamental representations of the quantum $\mathfrak{sl}(n)$.
Kuperberg generalized Temperley-Lieb algebras to
$\mathfrak{sl}(3)$ web spaces \cite{Kuperberg:spiders},
In section \ref{relation} we develop the quantum
$\mathfrak{sl}(n)$ representation theory by extending
the idea of webs in~\cite{MOY:Homfly}. Recently, Westbury found a web space for
spin representations of $\mathfrak{so}(7)$~\cite{math.QA/0601209}.
A precise and algebraic overview of the quantum $\mathfrak{sl}(n)$ representation theory can be found in ~\cite{Morrison}.
By expanding all crossings as in Figure~\ref{expansion}, Murakami, Ohtsuki and Yamada
found a regular isotopy invariant $[D]_n$~\cite{MOY:Homfly}.
In section~\ref{invariants}, we modify writhes suitably to define
an isotopy invariant $K_n(L,\mu)$, the $\mathfrak{sl}(n)$ HOMFLY polynomial,
where $\mu$ is a coloring of $L$ by a fundamental representation of the quantum $\mathfrak{sl}(n)$.
Using the quantum $\mathfrak{sl}(n)$ representation theory, we show
$K_n(L,\mu)$ can be computed by a linear expansion with respect to
the relations of web spaces in Theorem~\ref{onedim}.

\begin{figure}$$
\begin{pspicture}[shift=-.8](0,-1)(0,.7) \end{pspicture} \begin{pspicture}[shift=-.8](0,.8)(0,.8)
\begin{pspicture}[shift=-.8](-1,-1)(1,.7) \rput(.5,.5){\rnode{a1}{$$}} \rput(-.5,.5){\rnode{a2}{$$}}
\rput(-.5,-.5){\rnode{a3}{$$}} \rput(.5,-.5){\rnode{a4}{$$}} \rput(.1,.1){\rnode{b1}{$$}} \rput(-.1,.1){\rnode{b2}{$$}}
\rput(-.1,-.1){\rnode{b3}{$$}} \rput(.1,-.1){\rnode{b4}{$$}} \ncline{a1}{b1}\middlearrow  \ncline{a2}{b2}\middlearrow
\ncline{b4}{a4} \ncline{b1}{a3} \rput(0,-.8){\rnode{c4}{$L_+$}} \end{pspicture}\end{pspicture} \quad \quad
\begin{pspicture}[shift=-.8](0,.8)(0,.8) \begin{pspicture}[shift=-.8](-1,-1)(1,.7) \rput(.5,.5){\rnode{a1}{$$}} \rput(-.5,.5){\rnode{a2}{$$}}
\rput(-.5,-.5){\rnode{a3}{$$}} \rput(.5,-.5){\rnode{a4}{$$}} \rput(.1,.1){\rnode{b1}{$$}} \rput(-.1,.1){\rnode{b2}{$$}} \rput(-.1,-.1){\rnode{b3}{$$}}
\rput(.1,-.1){\rnode{b4}{$$}} \ncline{a2}{b2}\middlearrow \ncline{a1}{b1}\middlearrow \ncline{b2}{a4} \ncline{b3}{a3}
\rput(0,-.8){\rnode{c4}{$L_-$}} \end{pspicture}\end{pspicture}
\quad  \quad \begin{pspicture}[shift=-.8](0,.8)(0,.8) \begin{pspicture}[shift=-.8](-1,-1)(1,.7)
\rput(.5,.5){\rnode{a1}{$$}} \rput(-.5,.5){\rnode{a2}{$$}}
\rput(-.5,-.5){\rnode{a3}{$$}} \rput(.5,-.5){\rnode{a4}{$$}}
\rput(.1,.1){\rnode{b1}{$$}} \rput(-.1,.1){\rnode{b2}{$$}}
\rput(-.1,-.1){\rnode{b3}{$$}} \rput(.1,-.1){\rnode{b4}{$$}}
\nccurve[angleA=225,angleB=135]{a1}{a4}\middlearrow
\nccurve[angleA=315,angleB=45]{a2}{a3}\middlearrow
\rput(0,-.8){\rnode{c4}{$L_0$}}
\end{pspicture}\end{pspicture}$$
\caption{The skein triple $L_+, L_-$ and $L_0$.} \label{local}
\end{figure}
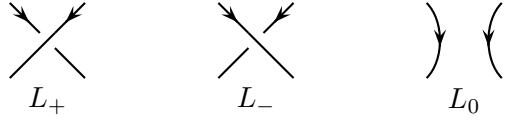

\begin{figure} $$
 \begin{pspicture}[shift=-.8](0,-.9)(0,.9) \end{pspicture}
 \begin{pspicture}[shift=-.8](0,.9)(0,.9)
\begin{pspicture}[shift=-.8](-.9,-.9)(.9,.9)\rput(.5,.5){\rnode{a1}{$$}} \rput(-.5,.5){\rnode{a2}{$$}}
\rput(-.5,-.5){\rnode{a3}{$$}} \rput(.5,-.5){\rnode{a4}{$$}} \rput(.1,.1){\rnode{b1}{$$}}
\rput(-.1,.1){\rnode{b2}{$$}} \rput(-.1,-.1){\rnode{b3}{$$}}
\rput(.1,-.1){\rnode{b4}{$$}} \ncline{a2}{b2}\middlearrow
\ncline{a1}{b1}\middlearrow \ncline{b4}{a4} \ncline{b1}{a3}
\rput(.55,.7){\rnode{c1}{$1$}} \rput(-.55,.7){\rnode{c1}{$1$}}
\end{pspicture}\end{pspicture} = q^{\frac{1}{2}}
 \begin{pspicture}[shift=-.8](0,.9)(0,.9)
\begin{pspicture}[shift=-.8](-.9,-.9)(.9,.9)  \rput(.5,.5){\rnode{a1}{$$}} \rput(-.5,.5){\rnode{a2}{$$}}
\rput(-.5,-.5){\rnode{a3}{$$}} \rput(.5,-.5){\rnode{a4}{$$}}
\nccurve[angleA=225,angleB=135]{a1}{a4}\middlearrow
\nccurve[angleA=315,angleB=45]{a2}{a3}\middlearrow
\rput(.55,.7){\rnode{c1}{$1$}} \rput(-.55,.7){\rnode{c1}{$1$}}
\end{pspicture}\end{pspicture} -
 \begin{pspicture}[shift=-.8](0,.9)(0,.9)
\begin{pspicture}[shift=-.8](-.9,-.9)(.9,.9)  \rput(.5,.5){\rnode{a1}{$$}} \rput(-.5,.5){\rnode{a2}{$$}}
\rput(-.5,-.5){\rnode{a3}{$$}} \rput(.5,-.5){\rnode{a4}{$$}} \rput(0,.3){\rnode{b1}{$$}}
\rput(0,-.3){\rnode{b2}{$$}} \ncline{a1}{b1}\middlearrow
\ncline{a2}{b1}\middlearrow \ncline{b1}{b2}\middlearrow
\ncline{b2}{a3}\middlearrow \ncline{b2}{a4}\middlearrow
\rput(.7,.7){\rnode{c1}{$1$}} \rput(-.7,.7){\rnode{c1}{$1$}}
\rput(-.7,-.7){\rnode{c1}{$1$}} \rput(.7,-.7){\rnode{c1}{$1$}}
\rput(.3,0){\rnode{c4}{$2$}}
\end{pspicture}\end{pspicture} $$ $$ \begin{pspicture}[shift=-.8](0,-.9)(0,.9) \end{pspicture}
 \begin{pspicture}[shift=-.8](0,.9)(0,.9)
\begin{pspicture}[shift=-.8](-.9,-.9)(.9,.9)\rput(.5,.5){\rnode{a1}{$$}} \rput(-.5,.5){\rnode{a2}{$$}}
\rput(-.5,-.5){\rnode{a3}{$$}} \rput(.5,-.5){\rnode{a4}{$$}} \rput(.1,.1){\rnode{b1}{$$}}
\rput(-.1,.1){\rnode{b2}{$$}} \rput(-.1,-.1){\rnode{b3}{$$}}
\rput(.1,-.1){\rnode{b4}{$$}} \ncline{a2}{b2}\middlearrow
\ncline{a1}{b1}\middlearrow \ncline{b2}{a4} \ncline{b3}{a3}
\rput(.55,.7){\rnode{c1}{$1$}} \rput(-.55,.7){\rnode{c1}{$1$}}
\end{pspicture}\end{pspicture} = q^{-\frac{1}{2}}
 \begin{pspicture}[shift=-.8](0,.9)(0,.9)
\begin{pspicture}[shift=-.8](-.9,-.9)(.9,.9) \rput(.5,.5){\rnode{a1}{$$}} \rput(-.5,.5){\rnode{a2}{$$}}
\rput(-.5,-.5){\rnode{a3}{$$}} \rput(.5,-.5){\rnode{a4}{$$}}
\nccurve[angleA=225,angleB=135]{a1}{a4}\middlearrow
\nccurve[angleA=315,angleB=45]{a2}{a3}\middlearrow
\rput(.55,.7){\rnode{c1}{$1$}} \rput(-.55,.7){\rnode{c1}{$1$}}
\end{pspicture}\end{pspicture} -
 \begin{pspicture}[shift=-.8](0,.9)(0,.9)
\begin{pspicture}[shift=-.8](-.9,-.9)(.9,.9)
\rput(.5,.5){\rnode{a1}{$$}} \rput(-.5,.5){\rnode{a2}{$$}}
\rput(-.5,-.5){\rnode{a3}{$$}} \rput(.5,-.5){\rnode{a4}{$$}}  \rput(0,.2){\rnode{b1}{$$}}
\rput(0,-.2){\rnode{b2}{$$}} \ncline{a1}{b1}\middlearrow
\ncline{a2}{b1}\middlearrow \ncline{b1}{b2}\middlearrow
\ncline{b2}{a3}\middlearrow \ncline{b2}{a4}\middlearrow
\rput(.7,.7){\rnode{c1}{$1$}} \rput(-.7,.7){\rnode{c1}{$1$}}
\rput(-.7,-.7){\rnode{c1}{$1$}} \rput(.7,-.7){\rnode{c1}{$1$}}
\rput(.3,0){\rnode{c4}{$2$}}
\end{pspicture}\end{pspicture}$$\caption{Expansions of crossings for $P_n(L)$.} \label{planar}
\end{figure}

For $n=2$, we can decorate $L$ by any other
irreducible representations $V_n$ using the highest-weight projection
$$f_n:V_1^{\otimes n} \to V_1^{\otimes n}$$
whose image is $V_n$ where $V_1$ is the vector representation of $\mathfrak{sl}(2)$.
This projection is called a \emph{Jones-Wenzl projector} \cite{Wenzl:Proj}.
It does exist for $n\ge 3$ and called a \emph{clasp}
\cite{Kuperberg:spiders}.
Using these clasps, Lickorish first found a quantum
$\mathfrak{sl}(2)$ invariants of $3$-manifolds \cite{Lickorish:su}.
Ohtsuki and Yamada generalized it for the quantum $\mathfrak
{sl}(3)$ \cite{OY:quantum} and Yokota did for the quantum
$\mathfrak{sl}(n)$ \cite{yokota:skeinforn}.
A benefit of using the quantum $\mathfrak{sl}(n)$ representation theory for link
invariants is that some nontrivial facts from the original work of
\cite{Jones:braid, RT:1, RT:2} do follow easily such as the
integrality~\cite{Le:integral}.

By decorating each component by $\mu$,
we can define the colored $\mathfrak{sl}(n)$
HOMFLY polynomial, denoted by $G_n(L,\mu)$ as follows.
 For a given colored link $L$ of $l$ components
say, $L_1, L_2, \ldots, L_l$, where each component $L_i$ is colored
by an irreducible representation $V_{a(i)_1\lambda_1 +
a(i)_2\lambda_2 + \ldots + a(i)_{n-1}\lambda_{n-1}}$ of
$\mathfrak{sl}(n)$ and $\lambda_1$, $\lambda_2$, $\ldots$, $\lambda_{n-1}$
are the fundamental weights of $\mathfrak{sl}(n)$. The coloring is
denoted by $\mu=(a(1)_1\lambda_1 + a(1)_2\lambda_2 + \ldots +
a(1)_{n-1}\lambda_{n-1}, a(2)_1\lambda_1 + a(2)_2\lambda_2 + \ldots
+ a(2)_{n-1}\lambda_{n-1}, \ldots, a(l)_1\lambda_1 + a(l)_2\lambda_2
+ \ldots + a(l)_{n-1}\lambda_{n-1})$. First we replace each
component $L_i$ by $a(i)_1+a(i)_2+\ldots +a(i)_{n-1}$ copies of
parallel lines and each $a(i)_j$ line is colored by the weight
$\lambda_j$. Then we put a clasp of weight $(a(i)_1\lambda_1 +
a(i)_2\lambda_2 + \ldots + a(i)_{n-1}\lambda_{n-1})$ for $L_i$. If
we assume the clasps are far away from crossings, we expand each
crossing as in Figure~\ref{expansion}, then
clasps \cite{yokota:skeinforn}. The value we obtained after removing all
faces by using the relations is \emph{the colored $\mathfrak{sl}(n)$
HOMFLY polynomial} $G_n(L,\mu)$ of $L$.

\section{The quantum $\mathfrak{sl}(n)$ representation theory}
\label{relation}

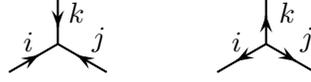
\begin{figure}
$$
 \begin{pspicture}[shift=-.6](0,-.4)(0,.8)\end{pspicture} \begin{pspicture}[shift=-.6](0,.6)(0,.6)
\begin{pspicture}[shift=-.6](-1,-.4)(1,.8) \rput(0,.6){\rnode{a1}{$$}}
\rput(-.65,-.375){\rnode{a2}{$$}} \rput(.65,-.375){\rnode{a3}{$$}}
\rput(0,0){\rnode{b1}{$$}}\rput[bl](4pt,8pt){$k$}
\rput[bl](-13pt,-3pt){$i$} \rput[bl](13pt,-3pt){$j$}
\ncline{a1}{b1}\middlearrow \ncline{a2}{b1}\middlearrow
\ncline{a3}{b1}\middlearrow
\end{pspicture}\end{pspicture}
\quad \quad \begin{pspicture}[shift=-.6](0,.6)(0,.6)
\begin{pspicture}[shift=-.6](-1,-.4)(1,.8)
\rput(0,.6){\rnode{a1}{$$}}
\rput(-.65,-.375){\rnode{a2}{$$}} \rput(.65,-.375){\rnode{a3}{$$}}
\rput(0,0){\rnode{b1}{$$}}
\rput[bl](5pt,8pt){$k$} \rput[bl](-13pt,-3pt){$i$}
\rput[bl](13pt,-3pt){$j$}
 \ncline{b1}{a1}\middlearrow
\ncline{b1}{a2}\middlearrow \ncline{b1}{a3}\middlearrow
\end{pspicture}\end{pspicture}
$$
\caption{Generators of the quantum $\mathfrak{sl}(n)$ web space.}
\label{generator}
\end{figure}

We refer to \cite{FultonHarris:gtm, Kuperberg:spiders, KRT:knot,
KR:factor} for the general representation theory. The webs are
generated by the two types of shapes in Figure~\ref{generator} where
$i+j+k=n$. At the second web in Figure~\ref{generator}, we can
change the directions of the edges to inward by using the duality,
$(V_{\lambda_i})^{*}\cong V_{\lambda_{n-i}}$, then the colors of the
resulting edges are $n-i, n-j$ and $n-k$. Then we say a vertex
has a {\it sign type} $+$ if the sum of the colors is $n$, $-$ if
the sum is $2n$. The generators in Figure~\ref{generator} have sign
type $+, -$ from left. Naturally we can assign a sign type for each
face. Using last two relations in
Figure~\ref{firstrel}, we only need to look at the face of a
\emph{valid sign type} in which $+, -$ are appearing alternatively.
One can easily see that none of faces with odd numbers of edges can have
a valid sign type.

A few computations of tensor products of fundamental representations
show the relations in Figure~\ref{firstrel} (up to scaling): the
quantum dimension of the
fundamental representation $V_{\lambda_i}$ is $\left[ \begin{matrix} n\\
i \end{matrix} \right]$, the dimension of the invariant space of
$V_{\lambda_i} \otimes V_{\lambda_i}^*$ is $1$ and the dimension of
the invariant space of $V_{\lambda_i} \otimes V_{\lambda_j}\otimes
V_{\lambda_k} \otimes V_{\lambda_i+\lambda_j+\lambda_k}^*$ is $1$.

\subsection{Rectangular relations}

To discuss rectangular relations, we first prove the following
lemma.

\begin{lem}
\begin{enumerate}
\item If $i\le j\le n-j-1$, then
$$\mathrm{dim}(\mathrm{inv}(V_{\lambda_i}\otimes V_{\lambda_j}
\otimes V_{\lambda_i}^{*} \otimes V_{\lambda_j}^{*}))=i+1.$$
\item
If $j\ge i\ge k\ge 1, n-j-1\ge i$ and $n-i-j-1\ge l\ge 1$, then
$$\mathrm{dim}(\mathrm{inv}(V_{\lambda_i}\otimes V_{\lambda_{j+l}}
\otimes V_{\lambda_{i+l}}^{*} \otimes V_{\lambda_{j}}^{*}))=i+1.$$
\end{enumerate}
\label{twotensor}
\end{lem}
\begin{proof}
If $i\le j\le n-j-1$, we obtain the following isomorphism by the
Clebsch-Gordan formula.

$$
V_{\lambda_i}\otimes V_{\lambda_j}\cong
V_{\lambda_{j-i}}\oplus V_{\lambda_{j-i+2}} \oplus
\cdots \oplus V_{\lambda_{j+i}}.$$

For irreducible representations $V, W$ of a simple Lie algebra, by a
simple application of Schur's lemma we find

$$\mathrm{dim}(\mathrm{inv}( V\otimes W^*)) = \left\{
\begin{array}{cl} 1 & ~~\mathrm{if}~~ V\cong W, \\
0 & ~~\mathrm{if}~~ V\not\cong W. \end{array}\right.
$$

These two facts imply that
$\mathrm{dim}(\mathrm{inv}(V_{\lambda_i}\otimes V_{\lambda_j}
\otimes V_{\lambda_i}^{*} \otimes V_{\lambda_j}^{*}))=i+1.$
Similarly one can prove the other.
\end{proof}
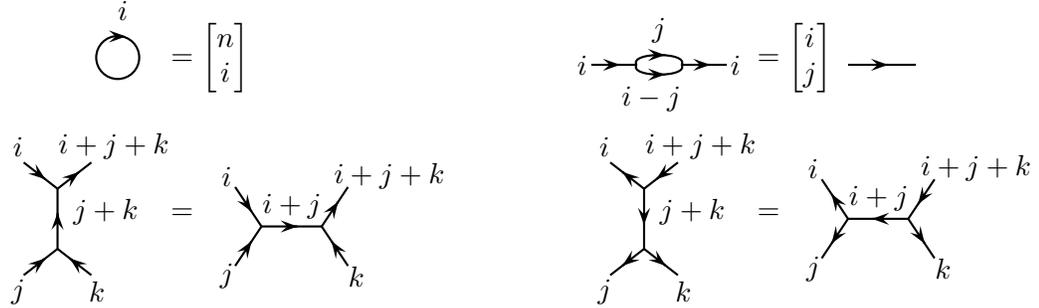
\begin{figure}
\begin{align*}
\begin{pspicture}[shift=-.3](-.6,-.4)(.6,.4) \pscircle(0,0){.3}
\psline[arrowscale=1.5]{->}(-.01,.28)(.11,.28)
\rput[bl](0pt,14pt){$i$}
\end{pspicture}
&=  \left[ \begin{matrix} n\\
i \end{matrix} \right] \quad \quad &
\begin{pspicture}[shift=-.6](-0,-.7)(0,.7)\end{pspicture}
\begin{pspicture}[shift=-.5](0,.6)(0,.6)
\begin{pspicture}[shift=-.5](-1.2,-.6)(1.2,.6)
\rput(-.9,0){\rnode{a1}{$$}}\rput(.9,0){\rnode{a4}{$$}}
\rput(-1.02,0){\rnode{b1}{$i$}}\rput(1.02,0){\rnode{b4}{$i$}}
\pnode(-.3,0){a2} \pnode(.3,0){a3} \ncline{a1}{a2}\middlearrow
\nccurve[angleA=120,angleB=60]{a2}{a3}\middlearrow
\rput(0pt,13pt){$j$}
\nccurve[angleA=-120,angleB=-60]{a2}{a3}\middlearrow
\rput(-3pt,-13pt){$i-j$} \ncline{a3}{a4}\middlearrow
\end{pspicture}\end{pspicture}
 &= \left[ \begin{matrix} i\\ j \end{matrix}
\right] \begin{pspicture}[shift=-.5](0,.6)(0,.6)
\begin{pspicture}[shift=-.5](-.6,-.6)(.6,.6) \rput(-.45,0){\rnode{a1}{$$}}
\rput(.45,0){\rnode{a2}{$$}} \ncline{a1}{a2}\middlearrow
\end{pspicture}\end{pspicture}  \\
\begin{pspicture}[shift=-1.1](-0,-1.2)(0,1.2)\end{pspicture}
\begin{pspicture}[shift=-1.1](0,1.2)(0,1.2)
\begin{pspicture}[shift=-1.1](-.8,-1.2)(1.4,1.2) \rput(-.45,.75){\rnode{c1}{$$}}
\rput(-.45,-.75){\rnode{c2}{$$}} \rput(.45,-.75){\rnode{c3}{$$}}
\rput(.45,.75){\rnode{c4}{$$}} \rput(-.525,.9625){\rnode{d1}{$i$}}
\rput(-.525,-.9625){\rnode{d2}{$j$}} \rput(.525,-.9625){\rnode{d3}{$k$}}
\rput(.75,.9625){\rnode{d4}{$i+j+k$}} \rput(0,.4){\rnode{b1}{$$}}
\rput(0,-.4){\rnode{b2}{$$}} \ncline{c1}{b1}\middlearrow
\ncline{b2}{b1}\middlearrow \rput[bl](6pt,-2pt){$j+k$}
\ncline{c2}{b2}\middlearrow \ncline{c3}{b2}\middlearrow
\ncline{b1}{c4}\middlearrow
\end{pspicture}\end{pspicture}
 &=  \begin{pspicture}[shift=-.9](0,1)(0,1)
\begin{pspicture}[shift=-.9](-1.2,-.9)(1.2,.9)
\rput(-.75,.525){\rnode{c1}{$$}} \rput(-.75,-.525){\rnode{c2}{$$}}
\rput(.75,-.525){\rnode{c3}{$$}} \rput(.75,.525){\rnode{c4}{$$}}
\rput(-.8625,.6625){\rnode{d1}{$i$}} \rput(-.8625,-.6625){\rnode{d2}{$j$}}
\rput(.8625,-.6625){\rnode{d3}{$k$}} \rput(1.3,.6725){\rnode{d4}{$i+j+k$}}
\rput(-.4,0){\rnode{b1}{$$}} \rput(.4,0){\rnode{b2}{$$}}
\ncline{c1}{b1}\middlearrow \ncline{b1}{b2}\middlearrow
\rput[bl](-11pt,3pt){$i+j$} \ncline{c2}{b1}\middlearrow
\ncline{c3}{b2}\middlearrow \ncline{b2}{c4}\middlearrow
\end{pspicture}\end{pspicture} \quad\quad &
\begin{pspicture}[shift=-1.1](-0,-1.2)(0,1.2)\end{pspicture}
\begin{pspicture}[shift=-1.1](0,1.2)(0,1.2)
\begin{pspicture}[shift=-1.1](-.8,-1.2)(1.4,1.2)  \rput(-.45,.75){\rnode{c1}{$$}}
\rput(-.45,-.75){\rnode{c2}{$$}} \rput(.45,-.75){\rnode{c3}{$$}}
\rput(.45,.75){\rnode{c4}{$$}} \rput(-.525,.9625){\rnode{d1}{$i$}}
\rput(-.525,-.9625){\rnode{d2}{$j$}} \rput(.525,-.9625){\rnode{d3}{$k$}}
\rput(.75,.9625){\rnode{d4}{$i+j+k$}} \rput(0,.4){\rnode{b1}{$$}}
\rput(0,-.4){\rnode{b2}{$$}}\ncline{b1}{c1}\middlearrow
\ncline{b1}{b2}\middlearrow \rput[bl](6pt,-2pt){$j+k$}
\ncline{b2}{c2}\middlearrow \ncline{b2}{c3}\middlearrow
\ncline{c4}{b1}\middlearrow
\end{pspicture}\end{pspicture}
 &=  \begin{pspicture}[shift=-.9](0,.9)(0,.9)
\begin{pspicture}[shift=-.9](-1.2,-.9)(1.2,.9)
\rput(-.75,.525){\rnode{c1}{$$}} \rput(-.75,-.525){\rnode{c2}{$$}}
\rput(.75,-.525){\rnode{c3}{$$}} \rput(.75,.525){\rnode{c4}{$$}}
\rput(-.8625,.6625){\rnode{d1}{$i$}} \rput(-.8625,-.6625){\rnode{d2}{$j$}}
\rput(.8625,-.6625){\rnode{d3}{$k$}} \rput(1.3,.6725){\rnode{d4}{$i+j+k$}}
\rput(-.4,0){\rnode{b1}{$$}} \rput(.4,0){\rnode{b2}{$$}}
\ncline{b1}{c1}\middlearrow \ncline{b2}{b1}\middlearrow
\rput[bl](-11pt,3pt){$i+j$} \ncline{b1}{c2}\middlearrow
\ncline{b2}{c3}\middlearrow \ncline{c4}{b2}\middlearrow
\end{pspicture}\end{pspicture}
\end{align*}
\caption{Relations obtained by the quantum dimension and scaling.}
\label{firstrel}
\end{figure}

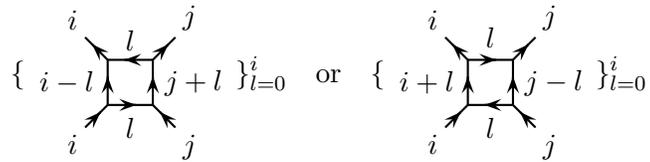
\begin{figure}
$$
\{ \begin{pspicture}[shift=-.8](0,-.9)(0,.9)
\end{pspicture} \begin{pspicture}[shift=-.8](0,.9)(0,.9)
\begin{pspicture}[shift=-.8](-1.4,-.9)(1.4,.9)\rput(.3,.3){\rnode{a1}{$$}}
\rput(0,0){\rnode{a0}{$$}} \rput(-.3,.3){\rnode{a2}{$$}}
\rput(-.3,-.3){\rnode{a3}{$$}} \rput(.3,-.3){\rnode{a4}{$$}}
\rput(.6,.6){\rnode{b1}{$$}} \rput(-.6,.6){\rnode{b2}{$$}}
\rput(-.6,-.6){\rnode{b3}{$$}} \rput(.6,-.6){\rnode{b4}{$$}}
\rput[bl](.7,.7){$j$} \rput[br](-.7,.7){$i$}
\rput[tr](-.7,-.7){$i$} \rput[tl](.7,-.7){$j$}
\rput(0pt,16pt){$l$} \rput(24pt,0pt){$j+l$} \rput(0pt,-17pt){$l$}
\rput(-24pt,0pt){$i-l$} \ncline{a1}{b1}\middlearrow
\ncline{a2}{b2}\middlearrow \ncline{b3}{a3}\middlearrow
\ncline{b4}{a4}\middlearrow \ncline{a1}{a2}\middlearrow
\ncline{a4}{a1}\middlearrow \ncline{a3}{a2}\middlearrow
\ncline{a3}{a4}\middlearrow
\end{pspicture}\end{pspicture}
\}_{l=0}^{i} \quad \mathrm{or} \quad \{
  \begin{pspicture}[shift=-.8](0,.9)(0,.9)
\begin{pspicture}[shift=-.8](-1.4,-.9)(1.4,.9)\rput(.3,.3){\rnode{a1}{$$}}
\rput(0,0){\rnode{a0}{$$}} \rput(-.3,.3){\rnode{a2}{$$}}
\rput(-.3,-.3){\rnode{a3}{$$}} \rput(.3,-.3){\rnode{a4}{$$}}
\rput(.6,.6){\rnode{b1}{$$}} \rput(-.6,.6){\rnode{b2}{$$}}
\rput(-.6,-.6){\rnode{b3}{$$}} \rput(.6,-.6){\rnode{b4}{$$}}
\rput[bl](.7,.7){$j$} \rput[br](-.7,.7){$i$}
\rput[tr](-.7,-.7){$i$} \rput[tl](.7,-.7){$j$}
\rput(0,.6){$l$} \rput(-24pt,0pt){$i+l$} \rput(0pt,-17pt){$l$}
\rput(24pt,0pt){$j-l$} \ncline{a1}{b1}\middlearrow
\ncline{a2}{b2}\middlearrow \ncline{b3}{a3}\middlearrow
\ncline{b4}{a4}\middlearrow \ncline{a2}{a1}\middlearrow
\ncline{a4}{a1}\middlearrow \ncline{a3}{a2}\middlearrow
\ncline{a4}{a3}\middlearrow
\end{pspicture}\end{pspicture}
\}_{l=0}^{i}
$$
\caption{Two basis webs with the boundary $V_{\lambda_i}\otimes
V_{\lambda_j}\otimes V_{\lambda_i}^* \otimes V_{\lambda_j}^*$. }
\label{recbasis}
\end{figure}

From Lemma~\ref{twotensor}, we know the number of basis webs that we
need for each expansion. For $n-i-1\ge j\ge i \ge 0$, we can have
two sets of basis webs and each has the same sign type as in Figure
\ref{recbasis}. Throughout the section we will use the basis in the
left hand side of Figure~\ref{recbasis}. There are only two possible
types of rectangular relations as in Lemma ~\ref{rectangular}, all
other can be taken care of by relations in Figure~\ref{firstrel}. The equation (\ref{recrel1})
in Lemma \ref{rectangular} was first appeared in \cite{MOY:Homfly} without
a proof.

\begin{lem}
For $n-i-1 \ge j\ge i\ge k\ge 0$ and $n-i-j-1\ge l\ge 0$, we find
\begin{eqnarray}\begin{pspicture}[shift=-1.1](0,-1.2)(0,1.2)
\end{pspicture}
 \begin{pspicture}[shift=-.9](0,1)(0,1)
\begin{pspicture}[shift=-1.1](-1.7,-1.2)(2.3,1.2) \rput(.8,.3){\rnode{a1}{$$}}
\rput(0,0){\rnode{a0}{$$}} \rput(-.8,.3){\rnode{a2}{$$}}
\rput(-.8,-.3){\rnode{a3}{$$}} \rput(.8,-.3){\rnode{a4}{$$}}
\rput(1.1,.8){\rnode{b1}{$$}} \rput(-1.1,.8){\rnode{b2}{$$}}
\rput(-1.1,-.8){\rnode{b3}{$$}} \rput(1.1,-.8){\rnode{b4}{$$}}
\rput(1.6,.9){\rnode{c1}{$j+l$}} \rput(-1.3,.9){\rnode{c2}{$i$}}
\rput(-1.3,-.9){\rnode{c3}{$j$}} \rput(1.6,-.9){\rnode{c4}{$i+l$}}
\rput(0pt,20pt){$j+k-i$} \rput(-40pt,0pt){$j+k$}
\rput(0pt,-18pt){$k$} \rput(49pt,0pt){$i-k+l$}
\ncline{a1}{b1}\middlearrow \ncline{a2}{b2}\middlearrow
\ncline{b3}{a3}\middlearrow \ncline{b4}{a4}\middlearrow
\ncline{a2}{a1}\middlearrow \ncline{a4}{a1}\middlearrow
\ncline{a3}{a2}\middlearrow \ncline{a4}{a3}\middlearrow
\end{pspicture}\end{pspicture}
&=& \sum _{m=0}^{i}\left[ \begin{matrix} l\\
k-m \end{matrix} \right]\hskip .4cm
 \begin{pspicture}[shift=-.9](0,1)(0,1)
\begin{pspicture}[shift=-1.1](-1.7,-1.2)(2,1.2)\rput(.8,.3){\rnode{a1}{$$}}
\rput(0,0){\rnode{a0}{$$}} \rput(-.8,.3){\rnode{a2}{$$}}
\rput(-.8,-.3){\rnode{a3}{$$}} \rput(.8,-.3){\rnode{a4}{$$}}
\rput(1.1,.8){\rnode{b1}{$$}} \rput(-1.1,.8){\rnode{b2}{$$}}
\rput(-1.1,-.8){\rnode{b3}{$$}} \rput(1.1,-.8){\rnode{b4}{$$}}
\rput(1.6,.9){\rnode{c1}{$j+l$}} \rput(-1.3,.9){\rnode{c2}{$i$}}
\rput(-1.3,-.9){\rnode{c3}{$j$}} \rput(1.6,-.9){\rnode{c4}{$i+l$}}
\rput(0pt,20pt){$m$} \rput(-43pt,0pt){$i-m$}
\rput(0pt,-18pt){$j-i+m$} \rput(51pt,0pt){$j+l+m$}
\ncline{a1}{b1}\middlearrow \ncline{a2}{b2}\middlearrow
\ncline{b3}{a3}\middlearrow \ncline{b4}{a4}\middlearrow
\ncline{a1}{a2}\middlearrow \ncline{a4}{a1}\middlearrow
\ncline{a3}{a2}\middlearrow \ncline{a3}{a4}\middlearrow
\end{pspicture}\end{pspicture} \label{recrel1}
\end{eqnarray}
For $n-j-1 \ge i\ge
j\ge k\ge 1$ and $n-i-j-1\ge l\ge 1$, we have
\begin{eqnarray} \begin{pspicture}[shift=-1.1](0,-1.2)(0,1.2)
\end{pspicture}
 \begin{pspicture}[shift=-.9](0,1)(0,1)
\begin{pspicture}[shift=-1.1](-1.7,-1.2)(2.3,1.2)\rput(.8,.3){\rnode{a1}{$$}}
\rput(0,0){\rnode{a0}{$$}} \rput(-.8,.3){\rnode{a2}{$$}}
\rput(-.8,-.3){\rnode{a3}{$$}} \rput(.8,-.3){\rnode{a4}{$$}}
\rput(1.1,.8){\rnode{b1}{$$}} \rput(-1.1,.8){\rnode{b2}{$$}}
\rput(-1.1,-.8){\rnode{b3}{$$}} \rput(1.1,-.8){\rnode{b4}{$$}}
\rput(1.6,.9){\rnode{c1}{$j+l$}} \rput(-1.3,.9){\rnode{c2}{$i$}}
\rput(-1.3,-.9){\rnode{c3}{$j$}} \rput(1.6,-.9){\rnode{c4}{$i+l$}}
\rput(0pt,20pt){$k$} \rput(-40pt,0pt){$i+k$}
\rput(0pt,-18pt){$i-j+k$} \rput(51pt,0pt){$j-k+l$}
\ncline{a1}{b1}\middlearrow \ncline{a2}{b2}\middlearrow
\ncline{b3}{a3}\middlearrow \ncline{b4}{a4}\middlearrow
\ncline{a2}{a1}\middlearrow \ncline{a4}{a1}\middlearrow
\ncline{a3}{a2}\middlearrow \ncline{a4}{a3}\middlearrow
\end{pspicture}\end{pspicture}
&=& \sum _{m=0}^{j}\left[ \begin{matrix} l\\
k-m \end{matrix} \right]\hskip .4cm
 \begin{pspicture}[shift=-.9](0,1)(0,1)
\begin{pspicture}[shift=-1.1](-1.7,-1.2)(2,1.2)\rput(.8,.3){\rnode{a1}{$$}}
\rput(0,0){\rnode{a0}{$$}} \rput(-.8,.3){\rnode{a2}{$$}}
\rput(-.8,-.3){\rnode{a3}{$$}} \rput(.8,-.3){\rnode{a4}{$$}}
\rput(1.1,.8){\rnode{b1}{$$}} \rput(-1.1,.8){\rnode{b2}{$$}}
\rput(-1.1,-.8){\rnode{b3}{$$}} \rput(1.1,-.8){\rnode{b4}{$$}}
\rput(1.6,.9){\rnode{c1}{$j+l$}} \rput(-1.3,.9){\rnode{c2}{$i$}}
\rput(-1.3,-.9){\rnode{c3}{$j$}} \rput(1.6,-.9){\rnode{c4}{$i+l$}}
\rput[a0](0pt,20pt){$i-j+m$} \rput[a0](-43pt,0pt){$j-m$}
\rput[a0](0pt,-18pt){$m$} \rput[a0](51pt,0pt){$i+l+m$}
\ncline{a1}{b1}\middlearrow \ncline{a2}{b2}\middlearrow
\ncline{b3}{a3}\middlearrow \ncline{b4}{a4}\middlearrow
\ncline{a1}{a2}\middlearrow \ncline{a4}{a1}\middlearrow
\ncline{a3}{a2}\middlearrow \ncline{a3}{a4}\middlearrow
\end{pspicture}\end{pspicture} \label{recrel2}
\end{eqnarray}
Just for these two equations, we use a
different convention of quantum integers that
$\left[ \begin{matrix} 0\\
0 \end{matrix} \right]=1$ but $\left[ \begin{matrix} 0\\
s \end{matrix} \right]=0$ if $s \neq 0$. \label{rectangular}
\end{lem}

\begin{proof} Let $a(k,m), b(k,m)$ be the coefficients in the righthand side of the
equation~(\ref{recrel1}) and (\ref{recrel2}). We induct on $(\mathrm{min}\{i,j\},k)$ in lexicographic
order. The key idea is to prove both equations simultaneously. If
$k=0$ and $j\ge i$, we find the equation in
Figure~\ref{induction1}. For the case $k=0$ and $i\ge j$, it is
identical except the weight on horizontal arrow is replaced by the
weight $i-j$ of the opposite direction. For $i=0$, we find the
equations in Figure~\ref{induction2}. One can do for the case
$j=0$ similarly.

Now we are set to proceed to the induction step. Let us look at the first
case $n-i-1 \ge j\ge i\ge k\ge 0$, $n-i-j-1\ge l\ge 0$. On the top
of each web in equation~(\ref{recrel1}), we can attach $i$ different $H$'s, given in
Figure~\ref{Hs} where $i\ge s\ge 1$. If $s=i$, we can easily get

\begin{figure}
$$\begin{pspicture}[shift=-.7](0,-.9)(0,.6)\end{pspicture}
\begin{pspicture}[shift=-.5](0,.6)(0,.6)
\begin{pspicture}[shift=-.7](-1.1,-.9)(1.2,.6) \rput(-.75,.525){\rnode{c1}{$$}}
\rput(-.75,-.525){\rnode{c2}{$$}} \rput(.75,-.525){\rnode{c3}{$$}}
\rput(.75,.525){\rnode{c4}{$$}} \rput(-.825,.7){\rnode{d1}{$i$}}
\rput(-.825,-.7){\rnode{d2}{$j$}} \rput(.825,-.7){\rnode{d3}{$i+l$}}
\rput(.825,.7){\rnode{d4}{$j+l$}} \rput(-.45,0){\rnode{b1}{$$}}
\rput(.45,0){\rnode{b2}{$$}} \ncline{b1}{c1}\middlearrow
\ncline{b1}{b2}\middlearrow \rput[b](0pt,3pt){$j-i$}
\ncline{c2}{b1}\middlearrow \ncline{b2}{c4}\middlearrow
\ncline{c3}{b2}\middlearrow
\end{pspicture}\end{pspicture}
=  \left[ \begin{matrix} l\\
0 \end{matrix} \right] \begin{pspicture}[shift=-.5](0,.6)(0,.6)
\begin{pspicture}[shift=-.7](-1.1,-.9)(1.2,.6) \rput(-.75,.525){\rnode{c1}{$$}}
\rput(-.75,-.525){\rnode{c2}{$$}} \rput(.75,-.525){\rnode{c3}{$$}}
\rput(.75,.525){\rnode{c4}{$$}} \rput(-.825,.7){\rnode{d1}{$i$}}
\rput(-.825,-.7){\rnode{d2}{$j$}} \rput(.825,-.7){\rnode{d3}{$i+l$}}
\rput(.825,.7){\rnode{d4}{$j+l$}} \rput(-.45,0){\rnode{b1}{$$}}
\rput(.45,0){\rnode{b2}{$$}} \ncline{b1}{c1}\middlearrow
\ncline{b1}{b2}\middlearrow \rput[b](0pt,3pt){$j-i$}
\ncline{c2}{b1}\middlearrow \ncline{b2}{c4}\middlearrow
\ncline{c3}{b2}\middlearrow
\end{pspicture}\end{pspicture}
$$
\caption{Initial induction step $k=0$ and $j\ge i$.}
\label{induction1}
\end{figure}

\begin{figure}
$$
\begin{pspicture}[shift=-.6](0,-.7)(0,.7)\end{pspicture}
\begin{pspicture}[shift=-.6](0,.7)(0,.7)
\begin{pspicture}[shift=-.6](-1,-.7)(1.4,.7)\rput(.25,.25){\rnode{a1}{$$}}
\rput(0,0){\rnode{a0}{$$}} \rput(-.25,.25){\rnode{a2}{$$}}
\rput(-.25,-.25){\rnode{a3}{$$}} \rput(.25,-.25){\rnode{a4}{$$}}
\rput(.5,.5){\rnode{b1}{$$}} \rput(-.5,-.5){\rnode{b3}{$$}}
\rput(.5,-.5){\rnode{b4}{$$}} \rput[bl](.6,.6){$j+l$}
\rput[tr](-.6,-.6){$j$} \rput[tl](.6,-.6){$l$}
\rput[r](-.45,0){$j+k$} \rput[t](0,-.45){$k$}
\rput[l](.5,0){$l-k$} \ncline{a1}{b1}\middlearrow
\ncline{b3}{a3}\middlearrow \ncline{b4}{a4}\middlearrow
\ncline{a4}{a1}\middlearrow
\nccurve[angleA=90,angleB=180]{a3}{a1}\middlearrow
\ncline{a4}{a3}\middlearrow
\end{pspicture}\end{pspicture}
= \begin{pspicture}[shift=-.6](0,.7)(0,.7)
\begin{pspicture}[shift=-.6](-1,-.7)(1.4,.7)\rput(.25,.25){\rnode{a1}{$$}}
\rput(-.4,.4){\rnode{a0}{$$}} \rput(-.1,.1){\rnode{a40}{$$}}
\rput(.2,-.2){\rnode{a43}{$$}} \rput(.5,.5){\rnode{b1}{$$}}
\rput(-.5,-.5){\rnode{b3}{$$}} \rput(.5,-.5){\rnode{a44}{$$}}
\rput[bl](.6,.6){$j+l$}
\rput[tr](-.6,-.6){$j$} \rput[tl](.6,-.6){$l$} \rput[a0](-3pt,-10pt){$k$}
\rput[l](.4,0){$l-k$}
\nccurve[angleA=70,angleB=180]{a0}{b1}\middlearrow
\nccurve[angleA=90,angleB=200]{b3}{a0}\middlearrow
\ncline{a44}{a43}\middlearrow \ncline{a40}{a0}\middlearrow
\nccurve[angleA=45,angleB=45]{a43}{a40}\middlearrow
\nccurve[angleA=225,angleB=225]{a43}{a40}\middlearrow
\end{pspicture}\end{pspicture}
= \left[ \begin{matrix} l\\
k \end{matrix} \right]
\begin{pspicture}[shift=-.6](0,.7)(0,.7)
\begin{pspicture}[shift=-.6](-.8,-.7)(1,.7)\rput(.25,.25){\rnode{a1}{$$}}
\rput(0,0){\rnode{a0}{$$}} \rput(.5,.5){\rnode{b1}{$$}}
\rput(-.5,-.5){\rnode{b3}{$$}} \rput(.5,-.5){\rnode{a44}{$$}}
\rput[bl](.6,.6){$j+l$}
\rput[tr](-.6,-.6){$j$} \rput[tl](.6,-.6){$l$}
\nccurve[angleA=90,angleB=225]{a0}{b1}\middlearrow
\nccurve[angleA=45,angleB=180]{b3}{a0}\middlearrow \ncline{a44}{a0}\middlearrow
\end{pspicture}\end{pspicture}
$$
\caption{Initial induction step $i=0$.} \label{induction2}
\end{figure}
\begin{figure}
$$ \begin{pspicture}[shift=-.3](0,-.2)(0,.6)\end{pspicture}
\begin{pspicture}[shift=-.3](0,.4)(0,.4)
\begin{pspicture}[shift=-.3](-1.3,-.2)(1.3,.6)\rput(.5,.25){\rnode{a1}{$$}}
\rput(0,0){\rnode{a0}{$$}} \rput(-.5,.25){\rnode{a2}{$$}}
\rput(-.5,-.1){\rnode{a3}{$$}} \rput(.5,-.1){\rnode{a4}{$$}}
\rput(1,.5){\rnode{b1}{$$}} \rput(-1,.5){\rnode{b2}{$$}}
\rput[bl](1.1,.6){$j+l+s$}
\rput[br](-1.1,.6){$i-s$} \rput(-1,-.5){$$}
\rput(1,-.5){\rnode{b4}{$$}} \rput[a0](0pt,14pt){$s$}
\rput[a0](-20pt,-4pt){$i$} \rput[a0](30pt,-4pt){$j+l$}
\ncline{a1}{b1}\middlearrow \ncline{a2}{b2}\middlearrow
\ncline{a2}{a1}\middlearrow \ncline{a4}{a1}\middlearrow
\ncline{a3}{a2}\middlearrow
\end{pspicture}\end{pspicture}
$$
\caption{H's attached on top of webs in Lemma~\ref{rectangular}.} \label{Hs}
\end{figure}

$$
\left[ \begin{matrix} j+k\\
i \end{matrix} \right] \left[ \begin{matrix} i+l\\
k \end{matrix} \right] = \sum_{m=0}^i a(i,m) \left[ \begin{matrix} j+l+m\\
m \end{matrix} \right] \left[ \begin{matrix} j\\
i-m \end{matrix} \right].
$$

From the case $s=1$, first we apply the bottom two relations in Figure~\ref{firstrel} at the upper rectangle of the web in the left hand side of the first equality in Figure~\ref{lhs1}
then the second relation in Figure~\ref{firstrel} to obtain the first equality. The second equality in Figure~\ref{lhs1} follows from the induction hypothesis of the equation~(\ref{recrel1}).

Next we look at each term in right hand
side of equation~(\ref{recrel1}) as in the first web on in Figure~\ref{rhs1}.
Now we can use the equation~(\ref{recrel2}) for the upper rectangle of the first web because of the induction hypothesis, where
the indices of the boundary are $i-1$, $i-m$, $(i-1)+(j+l+m-i+1)$ and $(i-m)+(j+l-i+m+1)$ from the northwest corner counter-clockwisely.
Since $k=1$ there are only two nonzero terms as in the right hand side of the first equality in Figure~\ref{rhs1} where $\alpha=[-i+j+l+m+1]$, $\beta=1$.
For the second web in the right hand side of the first equality in Figure~\ref{rhs1}, one can see that a similar step of relations, which was used in Figure~\ref{lhs1}
can be applied for the lower rectangle to get the next equality.

At last, by comparing coefficients of each basis
element, we get the following $i$ equations

$$
[-i+j+k+1]\left[ \begin{matrix} l+1\\
k-t \end{matrix} \right] = a(i,t) [-i+j+t+1]+a(i,t+1) [-i+j+l+t+2],
$$
where $i-1 \ge t \ge 0$. Since these $i+1$ equations are
independent, we plug in the answer to equations to check $a(i,m)=\left[\begin{matrix} l\\k-m\end{matrix}\right]$ is
correct. One can follow the proof for the second case, $n-j-1 \ge i\ge
j\ge k\ge 1$, $n-i-j-1\ge l\ge 1$.
\end{proof}

\begin{figure}
\begin{align*}\begin{pspicture}[shift=-1](0,-.8)(0,1.6)\end{pspicture}
&\begin{pspicture}[shift=-1](0,0)(0,0)\begin{pspicture}[shift=-1.1](-1,-.6)(2.6,1.6)
\psline(0,0)(0,.5)\psline[arrowscale=1.5]{->}(0,.15)(0,.35)
\psline(0,.5)(0,1)\psline[arrowscale=1.5]{->}(0,.65)(0,.85)
\psline(0,0)(1,0)\psline[arrowscale=1.5]{<-}(.4,0)(.6,0)
\psline(0,.5)(1,.5)\psline[arrowscale=1.5]{<-}(.6,.5)(.4,.5)
\psline(0,1)(1,1)\psline[arrowscale=1.5]{<-}(.6,1)(.4,1)
\psline(1,0)(1,.5)\psline[arrowscale=1.5]{->}(1,.15)(1,.35)
\psline(1,.5)(1,1)\psline[arrowscale=1.5]{->}(1,.65)(1,.85)
\psline(0,0)(-.25,-.25)\psline[arrowscale=1.5]{->}(-.2,-.2)(-.05,-.05)
\psline(1,0)(1.25,-.25)\psline[arrowscale=1.5]{->}(1.2,-.2)(1.05,-.05)
\psline(0,1)(-.25,1.25)\psline[arrowscale=1.5]{->}(-.05,1.05)(-.2,1.2)
\psline(1,1)(1.25,1.25)\psline[arrowscale=1.5]{->}(1.05,1.05)(1.2,1.2)
\rput[r](-.15,.25){$j+k$} \rput[r](-.15,.8){$i$}
\rput[l](1.15,.25){$i-k+l$} \rput[l](1.15,.8){$j+l$}
\rput[bl](.45,1.1){$1$}
\rput[bl](.45,-.4){$k$} \rput[br](-.0,1.3){$i-1$}
\rput[bl](1.3,1.3){$j+l+1$} \rput[tr](-.3,-.3){$j$}
\rput[tl](1.3,-.3){$i+l$}
\end{pspicture}\end{pspicture}
=\left[ \begin{matrix} j+k-i+1\\
1 \end{matrix} \right] \begin{pspicture}[shift=-1](0,-.3)(0,-.3)\begin{pspicture}[shift=-1.1](-1.1,-.6)(3.9,1.6)
\psline(0,0)(0,.75)\psline[arrowscale=1.5]{->}(0,.3)(0,.5)
\psline(2.25,0)(2.25,.75)\psline[arrowscale=1.5]{->}(2.25,.3)(2.25,.5)
\psline(0,0)(2.25,0)\psline[arrowscale=1.5]{->}(1.15,0)(.95,0)
\psline(0,.75)(2.25,.75)\psline[arrowscale=1.5]{->}(.95,.75)(1.15,.75)
\psline(0,0)(-.25,-.25)\psline[arrowscale=1.5]{->}(-.2,-.2)(-.05,-.05)
\psline(2.25,0)(2.5,-.25)\psline[arrowscale=1.5]{->}(2.45,-.2)(2.3,-.05)
\psline(0,.75)(-.25,1)\psline[arrowscale=1.5]{->}(-.05,.8)(-.2,.95)
\psline(2.25,.75)(2.5,1)\psline[arrowscale=1.5]{->}(2.3,.8)(2.45,.95)
\rput[b](-.5,1.1){$i-1$} \rput[b](3.2,1.1){$j+l+1$}
\rput[t](-.3,-.3){$j$} \rput[t](2.8,-.3){$i+l$}
\rput[r](-.15,.4){$j+k$} \rput[l](2.4,.4){$i-k+l$}
\rput[bl](1.1,-.4){$k$}\rput[bl](0,.9){$j+k-i+1$}
\end{pspicture}\end{pspicture}  \\ \begin{pspicture}[shift=-1](0,-.8)(0,1.6)\end{pspicture}
&=\left[ \begin{matrix} j+k-i+1\\
1 \end{matrix} \right] \sum_{m=0}^{i-1} \left[ \begin{matrix} l+1\\
k-m \end{matrix} \right] \begin{pspicture}[shift=-1](0,-.3)(0,-.3)\begin{pspicture}[shift=-1.1](-1.85,-.6)(3.4,1.6)
\psline(0,0)(0,.75)\psline[arrowscale=1.5]{->}(0,.3)(0,.5)
\psline(2.25,0)(2.25,.75)\psline[arrowscale=1.5]{->}(2.25,.3)(2.25,.5)
\psline(0,0)(2.25,0)\psline[arrowscale=1.5]{<-}(1.15,0)(.95,0)
\psline(0,.75)(2.25,.75)\psline[arrowscale=1.5]{<-}(.95,.75)(1.15,.75)
\psline(0,0)(-.25,-.25)\psline[arrowscale=1.5]{->}(-.2,-.2)(-.05,-.05)
\psline(2.25,0)(2.5,-.25)\psline[arrowscale=1.5]{->}(2.45,-.2)(2.3,-.05)
\psline(0,.75)(-.25,1)\psline[arrowscale=1.5]{->}(-.05,.8)(-.2,.95)
\psline(2.25,.75)(2.5,1)\psline[arrowscale=1.5]{->}(2.3,.8)(2.45,.95)
\rput[b](-.3,1.1){$i-1$} \rput[b](3.2,1.1){$j+l+1$}
\rput[t](-.3,-.3){$j$} \rput[t](2.8,-.3){$i+l$}
\rput[br](-.1,.2){$i-m-1$} \rput[bl](2.4,.2){$j+l+m+1$}
\rput[bl](1.1,.9){$m$}\rput[bl](-.1,-.5){$j-i+m+1$}
\end{pspicture}\end{pspicture}\end{align*}
\caption{The reduction on the left side of equation~(\ref{recrel1}).} \label{lhs1}
\end{figure}

\begin{figure}
\begin{align*}
&\begin{pspicture}[shift=-.1](-.6,-.6)(3.3,1.6)
\psline(0,0)(0,.5)\psline[arrowscale=1.5]{->}(0,.15)(0,.35)
\psline(0,.5)(0,1)\psline[arrowscale=1.5]{->}(0,.65)(0,.85)
\psline(0,0)(1.5,0)\psline[arrowscale=1.5]{->}(.65,0)(.85,0)
\psline(0,.5)(1.5,.5)\psline[arrowscale=1.5]{->}(.85,.5)(.65,.5)
\psline(0,1)(1.5,1)\psline[arrowscale=1.5]{<-}(.85,1)(.65,1)
\psline(1.5,0)(1.5,.5)\psline[arrowscale=1.5]{->}(1.5,.15)(1.5,.35)
\psline(1.5,.5)(1.5,1)\psline[arrowscale=1.5]{->}(1.5,.65)(1.5,.85)
\psline(0,0)(-.25,-.25)\psline[arrowscale=1.5]{->}(-.2,-.2)(-.05,-.05)
\psline(1.5,0)(1.75,-.25)\psline[arrowscale=1.5]{->}(1.7,-.2)(1.55,-.05)
\psline(0,1)(-.25,1.25)\psline[arrowscale=1.5]{->}(-.05,1.05)(-.2,1.2)
\psline(1.5,1)(1.75,1.25)\psline[arrowscale=1.5]{->}(1.55,1.05)(1.7,1.2)
\rput[br](-.2,.15){$i-m$} \rput[bl](-.3,.65){$i$}
\rput[l](1.65,.3){$j+l+m$} \rput[l](1.65,.8){$j+l$}
\rput[bl](-.1,-.45){$j-i+m$} \rput(.75,1.25){$1$}
\rput(.75,.75){$m$} \rput[br](-.3,1.3){$i-1$}
\rput[bl](1.8,1.3){$j+l+1$} \rput[tr](-.3,-.3){$j$}
\rput[tl](1.8,-.3){$i+l$}
\end{pspicture} \\
=\alpha \begin{pspicture}[shift=-.1](-1.4,-.6)(3.5,1.6)
\psline(0,0)(0,.5)\psline[arrowscale=1.5]{->}(0,.15)(0,.35)
\psline(1.5,0)(1.5,.5)\psline[arrowscale=1.5]{->}(1.5,.15)(1.5,.35)
\psline(0,0)(1.5,0)\psline[arrowscale=1.5]{->}(.65,0)(.85,0)
\psline(0,.5)(1.5,.5)\psline[arrowscale=1.5]{->}(.85,.5)(.65,.5)
\psline(0,0)(-.25,-.25)\psline[arrowscale=1.5]{->}(-.2,-.2)(-.05,-.05)
\psline(1.5,0)(1.75,-.25)\psline[arrowscale=1.5]{->}(1.7,-.2)(1.55,-.05)
\psline(0,.5)(-.25,.75)\psline[arrowscale=1.5]{->}(-.05,.55)(-.2,.7)
\psline(1.5,.5)(1.75,.75)\psline[arrowscale=1.5]{->}(1.55,.55)(1.7,.7)
\rput[br](-.3,.8){$i-1$} \rput[bl](1.7,.8){$j+l+1$}
\rput[tr](-.3,-.3){$j$} \rput[tl](1.7,-.3){$i+l$}
\rput[br](-.1,.15){$i-m$} \rput[bl](1.65,.15){$j+l+m$}
\rput[bl](.2,.65){$m-1$} \rput[bl](-.1,-.45){$j-i+m$}
\end{pspicture}
&+ \beta\begin{pspicture}[shift=-.3](-1.4,-.7)(3.7,1.6)
\psline(0,0)(0,.5)\psline[arrowscale=1.5]{->}(0,.15)(0,.35)
\psline(0,.5)(0,1)\psline[arrowscale=1.5]{->}(0,.65)(0,.85)
\psline(0,0)(1.5,0)\psline[arrowscale=1.5]{->}(.7,0)(.8,0)
\psline(0,.5)(1.5,.5)\psline[arrowscale=1.5]{->}(.7,.5)(.8,.5)
\psline(0,1)(1.5,1)\psline[arrowscale=1.5]{->}(.8,1)(.7,1)
\psline(1.5,0)(1.5,.5)\psline[arrowscale=1.5]{->}(1.5,.15)(1.5,.35)
\psline(1.5,.5)(1.5,1)\psline[arrowscale=1.5]{->}(1.5,.65)(1.5,.85)
\psline(0,0)(-.25,-.25)\psline[arrowscale=1.5]{->}(-.2,-.2)(-.05,-.05)
\psline(1.5,0)(1.75,-.25)\psline[arrowscale=1.5]{->}(1.7,-.2)(1.55,-.05)
\psline(0,1)(-.25,1.25)\psline[arrowscale=1.5]{->}(-.05,1.05)(-.2,1.2)
\psline(1.5,1)(1.75,1.25)\psline[arrowscale=1.5]{->}(1.55,1.05)(1.7,1.2)
\rput[br](-.1,.15){$i-m$} \rput[br](-.1,.7){$i-m-1$}
\rput[bl](1.6,.15){$j+l+m$} \rput[bl](1.6,.7){$j+l+m+1$}
\rput[bl](-.1,-.45){$j-i+m$} \rput[bl](.45,.6){$1$}
\rput[bl](.45,1.1){$m$} \rput[b](-.3,1.3){$i-1$}
\rput[bl](1.8,1.3){$j+l+1$} \rput[t](-.3,-.3){$j$}
\rput[tl](1.8,-.3){$i+l$}
\end{pspicture} \\
=\alpha \begin{pspicture}[shift=-.1](-1.4,-.6)(3.3,1.2)
\psline(0,0)(0,.5)\psline[arrowscale=1.5]{->}(0,.15)(0,.35)
\psline(1.5,0)(1.5,.5)\psline[arrowscale=1.5]{->}(1.5,.15)(1.5,.35)
\psline(0,0)(1.5,0)\psline[arrowscale=1.5]{->}(.7,0)(.8,0)
\psline(0,.5)(1.5,.5)\psline[arrowscale=1.5]{->}(.8,.5)(.7,.5)
\psline(0,0)(-.25,-.25)\psline[arrowscale=1.5]{->}(-.2,-.2)(-.05,-.05)
\psline(1.5,0)(1.75,-.25)\psline[arrowscale=1.5]{->}(1.7,-.2)(1.55,-.05)
\psline(0,.5)(-.25,.75)\psline[arrowscale=1.5]{->}(-.05,.55)(-.2,.7)
\psline(1.5,.5)(1.75,.75)\psline[arrowscale=1.5]{->}(1.55,.55)(1.7,.7)
\rput[b](-.3,.8){$i-1$} \rput[bl](1.8,.8){$j+l+1$}
\rput[t](-.3,-.3){$j$} \rput[tl](1.8,-.3){$i+l$}
\rput[br](-.1,.15){$i-m$} \rput[l](1.65,.3){$j+l+m$}
\rput[bl](0.2,.65){$m-1$} \rput[bl](-.1,-.45){$j-i+m$}
\end{pspicture} &+ \gamma \begin{pspicture}[shift=-.1](-.6,-.6)(3.5,1.2)
\psline(0,0)(0,.5)\psline[arrowscale=1.5]{->}(0,.15)(0,.35)
\psline(2.5,0)(2.5,.5)\psline[arrowscale=1.5]{->}(2.5,.15)(2.5,.35)
\psline(0,0)(2.5,0)\psline[arrowscale=1.5]{->}(1.15,0)(1.35,0)
\psline(0,.5)(2.5,.5)\psline[arrowscale=1.5]{->}(1.35,.5)(1.1,.5)
\psline(0,0)(-.25,-.25)\psline[arrowscale=1.5]{->}(-.2,-.2)(-.05,-.05)
\psline(2.5,0)(2.75,-.25)\psline[arrowscale=1.5]{->}(2.7,-.2)(2.55,-.05)
\psline(0,.5)(-.25,.75)\psline[arrowscale=1.5]{<-}(-.2,.7)(-.05,.55)
\psline(2.5,.5)(2.75,.75)\psline[arrowscale=1.5]{->}(2.55,.55)(2.7,.7)
\rput[b](-.3,.8){$i-1$} \rput[bl](2.8,.8){$j+l+1$}
\rput[t](-.3,-.3){$j$} \rput[tl](2.8,-.3){$i+l$}
\rput[bl](.1,.1){$i-m-1$} \rput[bl](2.7,.15){$j+l+m+1$}
\rput[bl](1.1,.6){$m$} \rput[bl](0,-.45){$j-i+m+1$}
\end{pspicture}
\end{align*}
\caption{The reduction on the right side of equation~(\ref{recrel1}).} \label{rhs1}
\end{figure}

\subsection{The quantum $\mathfrak{sl}(n)$ skein modules}

Skein modules were introduced independently by V. Turaev~\cite{Turaev:skein} and
J. Przytycki~\cite{PS:skein} as a $\mathbb{C}[A^{\pm1}]$-module associated to
a $3$-manifold $M$ generated by framed links inside $M$ with local relations
known as Kauffman relations. In the case of $M=\mathbb{S}^3$ this construction
reduces to the Jones polynomial and in the general case, the evaluation of the skein module
at the root of unity is known to fir with the Topological Quantum Field Theory constructed in
~\cite{BHMV}. It can be generalized for arbitrary Lie algebra, remind that the Jones polynomial
came from the representation category of the quantum $\mathfrak{sl}(2)$. By replacing
$M$ by $F \times [0,1]$, framed links by framed links with a color and Kauffman relation by relations given Figure~\ref{local} and Lemma~\ref{rectangular} one can obtain a $\mathbb{C}[A^{\pm1}]$-module. We call it
the \emph{quantum $\mathfrak{sl}(n)$ skein modules}. In this section, we concentrate
on the quantum $\mathfrak{sl}(n)$ skein modules and prove Theorem~\ref{onedim}.

The authors have been trying to find a complete relation of the
quantum $\mathfrak{sl}(n)$ representation theory, but we find
that the size of a polygon we have to
find a suitable relation is increasing as $n$ increases : there  a rectangular relation for the
quantum $\mathfrak{sl}(3)$, a hexagonal relation for the quantum
$\mathfrak{sl}(4)$ and an octagonal relation for the quantum $\mathfrak{sl}(6)$. In particular, a complete set of relations for
the quantum $\mathfrak{sl}(4)$ representation theory is conjectured
\cite{Dongseok:thesis}.

\begin{rem} For a given $n$ and sufficiently large $m$, $i. e., n\ll m$,
we conjecture that any $2n$ polygon of the sign type $(+,$ $-,$ $+,$
$\ldots,$ $-)$ can be expanded to a sum of webs of polygons of
smaller sizes and $2n$ polygons of the sign type $(-$, $+$, $-$,
$\ldots$, $+)$ by relations of the
quantum $\mathfrak{sl}(m)$ representation theory. \label{remcon}
\end{rem}

As mentioned before, a complete set of relations
for the representation theory of
the quantum $\mathfrak{sl}(n)$ which contain
our relations in Figure~\ref{firstrel}, Lemma~\ref{rectangular} and new relation called
'Kekul$\acute{e}$ relation', which was first found by the second author
for $\mathfrak{sl}(4)$~\cite{Dongseok:thesis},
is conjectured in~\cite{Morrison}. Furthermore, it was also shown that Remark~\ref{remcon} is false~\cite{Morrison}.
Without using extra relations found in~\cite{Morrison}, we can prove the following theorem.

\begin{thm} The quantum $\mathfrak{sl}(n)$ skein module of the plane or sphere
has dimension $1$. \label{onedim}
\end{thm}

\begin{proof}
If we look at these webs without decorations, they are directed,
weighted, trivalent and planar graphs. We will consider these graphs
on $\mathbb{S}^2$ instead of $\mathbb{R}^2$ and assume that
all webs are without boundary for the rest of proof, $i.e.$,
no vertex of valence $1$.

We will claim that the dimension of the web space
without boundary is $1$. Suppose it is not true, then there
exists a web for which we can not take a value in $\mathbb{C} [q^{\frac12},
q^{-\frac12}]$ by repeatedly applying relations for the representation theory of
the quantum $\mathfrak{sl}(n)$, say $D$. We assume that $D$ has the smallest number
of faces among all counterexamples. Since a face of size less than
$4$ can be removed by using relations found in
Figure~\ref{firstrel}, we assume that all faces in $D$ have at
least $4$ edges. Since it can be shown easily that there are finitely
many trivalent graphs with a fixed number of faces whose sizes are bigger
than $3$, we further assume that $D$ has the maximum
number of rectangular faces among counterexamples of the minimal number of faces. Let $V, E$ and $F$ be
the number of vertices, edges and faces in $D$ on
$\mathbb{S}^2$, respectively. Let $F_i$ be the number of faces with
$i$ edges. We can easily find the following equations:

\begin{eqnarray}
3V&=&2E, \label{graphformula1}\\
2&=&V-E+F, \label{graphformula2}\\
2E&=& \sum_{4\le i} iF_i. \label{graphformula3}
\end{eqnarray}

The outline of the proof is as follow. First, we prove the existence of rectangular face in $D$ in Lemma~\ref{existrectangle}.
Then we divide cases depending on the neighborhood of the rectangular faces.
Lemma~\ref{a3lem1} will show that if a rectangular face is adjacent to a pentagon, then the adjacent pentagon is unique, we call it an
\emph{isolated} rectangular face.
If a rectangular face is not isolated, we call it a \emph{non-isolated} rectangular face.
For non-isolated rectangular faces, we further divide them as follow; if the sizes of all adjacent faces are bigger than $6$, then
we call it a \emph{non-moveable} rectangular face. If a rectangular face is non-isolated and it is adjacent to a hexagon, then
we call it a \emph{moveable} rectangular face. Then we observe the neighborhood of rectangular faces, i) isolated rectangular faces in
Lemma~\ref{neighpenta}, (ii) moveable rectangular faces in Lemma~\ref{isolatedlem}. By
looking at all polygons whose sizes are bigger than $6$, we find a contradictory inequality to prove the theorem.

\begin{lem}
Let $D$ be a counterexample with the hypothesis,
the minimality of the number of faces and the maximality of the number of rectangular faces among counterexamples. Then,
there can not be two adjacent rectangles of a valid sign type in $D$.
\label{noadjacentrectangles}
\end{lem}
\begin{proof}
Suppose not, then there exist two adjacent rectangles of a valid sign type in $D$.
We apply an equation in Lemma~\ref{rectangular} to change one of adjacent rectangles to a linear combination of webs with the rectangular
faces with the opposite sign type. Then, two opposite edges of the other rectangle of adjacent rectangles have
the same signs at each ends. By applying the last two relations in Figure \ref{firstrel},
a web with two adjacent rectangles are now changed to a linear combination of webs with
a hexagon and a bigon as depicted in Figure~\ref{tworec}, consequently, just a hexagon. These processes result that the original web $D$ is a linear
combination of webs of less numbers of faces, however, since $D$ is a counterexample, at least one of webs in this linear combination
has to be a counterexample. But this can not be happened
because of the minimality of the number of faces in $D$.
\end{proof}

Let $I$ be the number of non-isolated
rectangular faces in $D$. The following lemma shows that
each pentagon is adjacent to a unique isolated rectangular face, $i.e.$, $F_4-I=F_5$.

\begin{lem}
Let $D$ be a counterexample with the hypothesis,
the minimality of the number of faces and the maximality of the number of rectangular faces among counterexamples.
If there exists a pentagonal face in $D$, the pentagonal face is adjacent to a unique
rectangular face. \label{a3lem1}
\end{lem}

\begin{proof}
Suppose the pentagon is not
adjacent to a rectangle, $i.e.,$ the sizes of all adjacent faces of the
pentagon are bigger than $4$. But by applying one of the last two relations
as shown in Figure \ref{firstrel} to the pentagonal face, it
can be changed to a rectangle without changing the total number of
faces but increasing the number of the rectangles by $1$ in $D$.
This contradicts the maximality of the number of the rectangular
faces in $D$. Therefore, it must be adjacent to a rectangle. To show the uniqueness of the adjacent
rectangle we consider the pentagon relations. There are five
possible shapes by rotating the web in
Figure~\ref{penta} by $\frac{2\pi}{5}$. By applying the
equations in Lemma \ref{rectangular}, one can see that any two out
of these five shapes are related as the equation as illustrated in Figure
~\ref{pentarel}, which is called a \emph{pentagon relation}. If two
rectangles are adjacent to a pentagon, we use a pentagon relation to make
a linear combinations of webs with two
adjacent rectangles and webs of one less number of faces.
Since $D$ is a counterexample, one of these webs in the linear combination
has to be a counterexample. However, a web with two
adjacent rectangles can not be a counterexample as described before.
A web with one less number of faces than that of $D$ can not be a counterexample
neither because of the minimality of the number of faces in $D$.
Therefore, the adjacent rectangle of the pentagon has to be unique.
\end{proof}

The rectangular faces play the key role in the proof as mentioned before. We will show the existence of a rectangular face
in the following lemma.

\begin{lem}
Let $D$ be a counterexample with the hypothesis,
the minimality of the number of faces and the maximality of the number of rectangular faces among counterexamples.
Then, there exists a rectangular face in $D$. \label{existrectangle}
\end{lem}
\begin{proof}
By Lemma~\ref{a3lem1},
if there does not exist a rectangular face in $D$, then the size of all faces are bigger than $5$,
$i. e.,$ $2E\ge 6F$.
By combining with equation~(\ref{graphformula1}) and~(\ref{graphformula2}), we easily find a contradiction as follows,
$$2=V-E+F \le \frac23 E- E + \frac13E =0.$$
\end{proof}

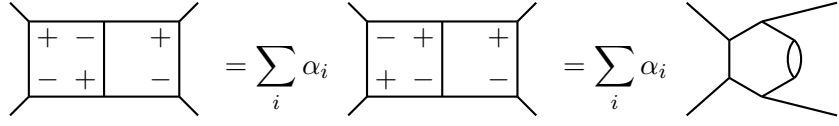
\begin{figure}
$$
\begin{pspicture}[shift=-.65](-.5,-.35)(2.5,1.35)
\psline(0,0)(0,1) \psline(2,0)(2,1)
\psline(0,0)(2,0) \psline(0,1)(2,1)
\psline(0,0)(-.25,-.25) \psline(2,0)(2.25,-.25)
\psline(0,1)(-.25,1.25) \psline(2,1)(2.25,1.25)
\psline(1,0)(1,1)
\rput[b](.25,.65){$+$} \rput[b](.75,.65){$-$}
\rput[b](.25,.1){$-$} \rput[b](.75,.1){$+$}
\rput[b](1.75,.65){$+$}  \rput[b](1.75,.1){$-$}
\end{pspicture}= \sum_i \alpha_i \begin{pspicture}[shift=-.65](-.5,-.35)(2.5,1.35)
\psline(0,0)(0,1) \psline(2,0)(2,1)
\psline(0,0)(2,0) \psline(0,1)(2,1)
\psline(0,0)(-.25,-.25) \psline(2,0)(2.25,-.25)
\psline(0,1)(-.25,1.25) \psline(2,1)(2.25,1.25)
\psline(1,0)(1,1)
\rput[b](.25,.65){$-$} \rput[b](.75,.65){$+$}
\rput[b](.25,.1){$+$} \rput[b](.75,.1){$-$}
\rput[b](1.75,.65){$+$}  \rput[b](1.75,.1){$-$}
\end{pspicture} = \sum_i \alpha_i \begin{pspicture}[shift=-.8](-1.25,-1)(1.1,1)
\rput(.433,.25){\rnode{a1}{$$}}
\rput(0,.5){\rnode{a2}{$$}}
\rput(-.433,.25){\rnode{a3}{$$}}
\rput(-.433,-.25){\rnode{a4}{$$}}
\rput(0,-.5){\rnode{a5}{$$}}
\rput(.433,-.25){\rnode{a6}{$$}}
\rput(1,.75){\rnode{b1}{$$}}
\rput(-1,.75){\rnode{b3}{$$}}
\rput(-1,-.75){\rnode{b4}{$$}}
\rput(1,-.75){\rnode{b6}{$$}}
\ncline{b1}{a2} \ncline{a3}{b3}
\ncline{a4}{b4} \ncline{a5}{b6} \ncline{a1}{a2} \ncline{a3}{a2}
\ncline{a3}{a4} \ncline{a5}{a4} \ncline{a5}{a6}
\nccurve[angleA=-45,angleB=45]{a1}{a6}
\nccurve[angleA=-135,angleB=135]{a1}{a6}
\end{pspicture}$$
\caption{Two adjacent rectangles with a valid sign type and its linear expansion of webs with less number of
faces.} \label{tworec}
\end{figure}

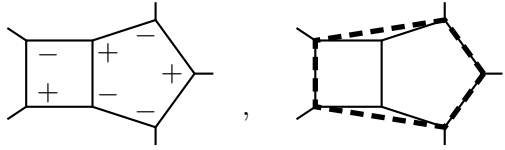
\begin{figure}
$$
\begin{pspicture}[shift=-.4](-1.9,-.95)(1,.95)
\psline(.75,0)(.2318,.7133)(-.6068,.4409)(-.6068,-.4409)(.2318,-.7133)(.75,0)
\psline(-.6068,.4409)(-1.4885,.4409)(-1.4885,-.4409)(-.6068,-.4409)
\qline(-1.4885,.4409)(-1.7385,.6068)
\qline(-1.4885,-.4409)(-1.7385,-.6068)
\qline(.2318,.7133)(.2318,.95)
\qline(.2318,-.7133)(.2318,-.95) \qline(.75,0)(1,0)
\rput(.45,0){$+$} \rput(.1,.5){$-$} \rput(.1,-.5){$-$}
 \rput(-.4,.28){$+$} \rput(-.4,-.28){$-$}
 \rput(-1.2,.25){$-$} \rput(-1.2,-.25){$+$}
\end{pspicture} \quad , \quad \begin{pspicture}[shift=-.4](-1.9,-.95)(1,.95)
\psline(.75,0)(.2318,.7133)(-.6068,.4409)(-.6068,-.4409)(.2318,-.7133)(.75,0)
\psline(-.6068,.4409)(-1.4885,.4409)(-1.4885,-.4409)(-.6068,-.4409)
\qline(-1.4885,.4409)(-1.7385,.6068)
\qline(-1.4885,-.4409)(-1.7385,-.6068)
\qline(.2318,.7133)(.2318,.95)
\qline(.2318,-.7133)(.2318,-.95) \qline(.75,0)(1,0)
\psline[linewidth=2pt,linestyle=dashed](.75,0)(.2318,.7133)(-1.4885,.4409)(-1.4885,-.4409)(.2318,-.7133)(.75,0)
\end{pspicture}$$
\caption{A valid sign type for an isolated rectangle and its adjacent pentagon
 and their corresponding new pentagon presented by dashed
lines considered in Lemma~\ref{neighpenta}.} \label{penta}
\end{figure}

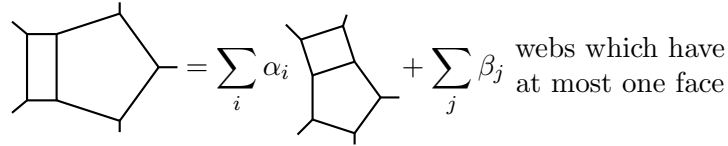
\begin{figure}
$$
\begin{pspicture}[shift=-.8](-1.3,-.9)(1,.9)
\psline(.75,0)(.2318,.7133)(-.6068,.4409)(-.6068,-.4409)(.2318,-.7133)(.75,0)
\psline(-.6068,.4409)(-1,.4409)(-1,-.4409)(-.6068,-.4409)
\qline(-1,.4409)(-1.2,.6068)
\qline(-1,-.4409)(-1.2,-.6068)
\qline(.2318,.7133)(.2318,.8633)
\qline(.2318,-.7133)(.2318,-.8633) \qline(.75,0)(1,0)
\end{pspicture}
= \sum_i\alpha_i\begin{pspicture}[shift=-.9](-.7,-.6)(.7,1.1)
\psline(.5,0)(.1545,.4756)(-.4045,.2939)(-.4045,-.2939)(.1545,-.4756)(.5,0)
\psline(-.4045,.2939)(-.559,.7699)(0,.9511)(.1545,.4756)
\qline(-.559,.7699)(-.759,.9699)
\qline(-.4045,-.2939)(-.6045,-.4939) \qline(0,.9511)(.05,1.1011)
\qline(.1545,-.4756)(.1545,-.6256) \qline(.5,0)(.75,0)
\end{pspicture} + \sum_j \beta_j ~~~\begin{matrix}\mathrm{webs}~\mathrm{which}~
\mathrm{have} \\ \mathrm{at}~
\mathrm{most}~ \mathrm{one}~
\mathrm{face}\end{matrix}
$$
\caption{A pentagonal relation.} \label{pentarel}
\end{figure}

\begin{figure}
$$
\begin{pspicture}[shift=-.6](-1.25,-.7)(1.1,.7)
\rput(.433,.25){\rnode{a1}{$$}}
\rput(0,.5){\rnode{a2}{$$}}
\rput(-.433,.25){\rnode{a3}{$$}}
\rput(-.433,-.25){\rnode{a4}{$$}}
\rput(0,-.5){\rnode{a5}{$$}}
\rput(.433,-.25){\rnode{a6}{$$}}
\rput(.65,.375){\rnode{b1}{$$}}
\rput(0,.75){\rnode{b2}{$$}}
\rput(-.65,.375){\rnode{b3}{$$}}
\rput(-.65,-.375){\rnode{b4}{$$}}
\rput(0,-.75){\rnode{b5}{$$}}
\rput(.65,-.375){\rnode{b6}{$$}}
\rput(.866,.25){\rnode{c1}{$$}}
\rput(.866,-.25){\rnode{c2}{$$}}
\rput(1.0427,.4267){\rnode{d1}{$$}}
\rput(1.0427,-.4267){\rnode{d2}{$$}}
\rput(0,0){\rnode{e7}{$$}} \ncline{b2}{a2} \ncline{a3}{b3}
\ncline{a4}{b4} \ncline{a5}{b5} \ncline{a1}{a2} \ncline{a3}{a2}
\ncline{a3}{a4} \ncline{a5}{a4} \ncline{a5}{a6} \ncline{a1}{a6}
\ncline{a1}{c1} \ncline{c2}{a6} \ncline{c2}{c1} \ncline{c1}{d1}
\ncline{c2}{d2}
\end{pspicture}
= \sum_i \alpha_i\hskip .5cm \begin{pspicture}[shift=-.6](-.9,-.7)(1.25,.7)
\rput(.433,.25){\rnode{a1}{$$}}
\rput(0,.5){\rnode{a2}{$$}}
\rput(-.433,.25){\rnode{a3}{$$}}
\rput(-.433,-.25){\rnode{a4}{$$}}
\rput(0,-.5){\rnode{a5}{$$}}
\rput(.433,-.25){\rnode{a6}{$$}}
\rput(.65,.375){\rnode{b1}{$$}}
\rput(0,.75){\rnode{b2}{$$}}
\rput(-.65,.375){\rnode{b3}{$$}}
\rput(-.65,-.375){\rnode{b4}{$$}}
\rput(0,-.75){\rnode{b5}{$$}}
\rput(.65,-.375){\rnode{b6}{$$}}
\rput(-.866,.25){\rnode{c1}{$$}}
\rput(-.866,-.25){\rnode{c2}{$$}}
\rput(-1.0427,.4267){\rnode{d1}{$$}}
\rput(-1.0427,-.4267){\rnode{d2}{$$}}
\rput(0,0){\rnode{e7}{$$}} \ncline{b2}{a2} \ncline{a1}{b1}
\ncline{a6}{b6} \ncline{a5}{b5} \ncline{a3}{a2} \ncline{a1}{a2}
\ncline{a4}{a3} \ncline{a5}{a6} \ncline{a5}{a4} \ncline{a1}{a6}
\ncline{c1}{a3} \ncline{a4}{c2} \ncline{c1}{c2} \ncline{c1}{d1}
\ncline{c2}{d2}
\end{pspicture} + \sum_j \beta_j ~~~\begin{matrix}\mathrm{webs}~\mathrm{which}~
\mathrm{have} \\ \mathrm{at}~
\mathrm{most}~ \mathrm{one}~
\mathrm{face}\end{matrix}
$$
\caption{A general swapping move for a moveable rectangular face.}\label{hexa3}
\end{figure}

To proceed the proof of Theorem~\ref{onedim}, we introduce a way modifying the web which has a moveable
rectangular face, recall that a moveable
rectangular face is a rectangular face adjacent to a hexagon. If a rectangular face is
moveable, the relations in Lemma~\ref{rectangular} allow
us to change the web $D$ to a linear combination of webs with a rectangular
face of an opposite sign type and webs of less number of faces.
For the webs with a rectangular
face of an opposite sign type, one can see that two edges
in the hexagonal face have the same sign at the ends.
By applying relations shown in Figure~\ref{firstrel}, the positions of
the rectangle and the hexagon are now interchanged as illustrated in
Figure~\ref{hexa3} which is called a \emph{swapping move}. Because of
the hypothesis of $D$, one of these new webs with adjacent rectangle and hexagon
has to be a counterexample, say $D'$. Slightly abusing notation, we will say that the original
counterexample $D$ can be changed to $D'$ by a swapping move.

In the following lemma, we look at the neighborhood of each isolated
rectangular face. Since every isolated rectangular face is adjacent to
a pentagon,
we will treat these two adjacent polygons as a new pentagon as depicted in the right side of Figure~\ref{penta}.

\begin{lem}
Let $D$ be a counterexample with the hypothesis,
the minimality of the number of faces and the maximality of the number of rectangular faces among counterexamples.
Suppose $D$ has an isolated rectangular face. Then,
the size of all adjacent polygons of a new pentagon drawn for an isolated
rectangular face and its adjacent pentagon as in the right side of Figure~\ref{penta} must be bigger than $7$.
\label{neighpenta}
\end{lem}

\begin{proof}
First we already know that these polygons adjacent to the pentagon can not be a rectangle nor
a pentagon. Suppose one of them is either a hexagon or a heptagon.
Either case, we can rotate the pentagon using the pentagon relations as shown in Figure~\ref{pentarel} such
that the rectangular face is adjacent to the hexagon or the heptagon (since the rectangle must have a
valid sign type, otherwise, it can be removed resulting a contradiction of the
minimality of the number of faces in $D$). But, we can change the heptagon to a hexagon using the relations
as illustrated in Figure~\ref{firstrel}. Then by a swapping move as depicted in
Figure~\ref{hexa3}, the rectangle in $D$ can be separated from the pentagon. Then
every web in the linear combinations has either one less number of faces or
a pentagon which is not adjacent to a rectangle, thus this pentagon can be changed to a rectangle using the relations
shown in Figure~\ref{firstrel}, consequently, it increases the number of rectangular faces. Since $D$ is a counterexample, at least one webs in the linear combinations
has to be a counterexample but neither cases is possible because of the
hypothesis of $D$.
\end{proof}

Now we proceed to the next lemma, we look at the neighborhood of a moveable
rectangular face. One can easily see that all adjacent hexagonal faces
must have a valid sign type (otherwise, they become two adjacent
rectangular faces using the relations
shown in Figure~\ref{firstrel}, which is impossible).
\begin{lem}
Let $D$ be a counterexample with the hypothesis,
the minimality of the number of faces and the maximality of the number of rectangular faces.
The only possible neighborhood of a moveable
rectangular face must be either one of webs in Figure~\ref{hexa1} or a circular web obtained from the right one in Figure~\ref{hexa1}
by attaching the rightmost and leftmost hexagons.
 \label{isolatedlem}
\end{lem}
\begin{proof}

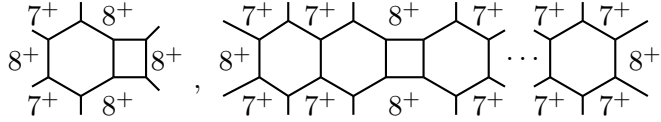
\begin{figure}
$$
\begin{pspicture}[shift=-.5](-1.25,-.9)(1.5,.9)
\rput(.433,.25){\rnode{a1}{$$}}
\rput(0,.5){\rnode{a2}{$$}}
\rput(-.433,.25){\rnode{a3}{$$}}
\rput(-.433,-.25){\rnode{a4}{$$}}
\rput(0,-.5){\rnode{a5}{$$}}
\rput(.433,-.25){\rnode{a6}{$$}}
\rput(.65,.375){\rnode{b1}{$$}}
\rput(0,.75){\rnode{b2}{$$}}
\rput(-.65,.375){\rnode{b3}{$$}}
\rput(-.65,-.375){\rnode{b4}{$$}}
\rput(0,-.75){\rnode{b5}{$$}}
\rput(.65,-.375){\rnode{b6}{$$}}
\rput(.866,.25){\rnode{c1}{$$}}
\rput(.866,-.25){\rnode{c2}{$$}}
\rput(1.0427,.4267){\rnode{d1}{$$}}
\rput(1.0427,-.4267){\rnode{d2}{$$}}
\rput(-.5,.6){\rnode{e1}{$7^+$}}
\rput(.5,.6){\rnode{e2}{$8^+$}}
\rput(-.5,-.6){\rnode{e3}{$7^+$}}
\rput(.5,-.6){\rnode{e4}{$8^+$}}
\rput(-.75,0){\rnode{e5}{$8^+$}}
\rput(1.15,0){\rnode{e4}{$8^+$}}
\rput(0,0){\rnode{e7}{$$}} \ncline{b2}{a2} \ncline{a3}{b3}
\ncline{a4}{b4} \ncline{a5}{b5} \ncline{a1}{a2} \ncline{a3}{a2}
\ncline{a3}{a4} \ncline{a5}{a4} \ncline{a5}{a6} \ncline{a1}{a6}
\ncline{a1}{c1} \ncline{c2}{a6} \ncline{c2}{c1} \ncline{c1}{d1}
\ncline{c2}{d2}
\end{pspicture},
 \quad
\begin{pspicture}[shift=-.5](-2.25,-.9)(3.35,.9) \rput(3.2,0){\rnode{e1}{$8^+$}}
\rput(-2.25,0){\rnode{e2}{$8^+$}}
\rput(-1.952,.6){\rnode{e3}{$7^+$}}
\rput(-1.116,.6){\rnode{e4}{$7^+$}} \rput(0,.6){\rnode{e5}{$8^+$}}
\rput(1.116,.6){\rnode{e6}{$7^+$}}
\rput(1.933,.6){\rnode{e8}{$7^+$}}
\rput(2.799,.6){\rnode{e9}{$7^+$}}
\rput(1.57,0){\rnode{f2}{$\ldots$}}
\rput(-1.952,-.6){\rnode{f3}{$7^+$}}
\rput(-1.116,-.6){\rnode{f4}{$7^+$}}
\rput(0,-.6){\rnode{f5}{$8^+$}}
\rput(1.116,-.6){\rnode{f6}{$7^+$}}
\rput(1.933,-.6){\rnode{f8}{$7^+$}}
\rput(2.799,-.6){\rnode{f9}{$7^+$}}
\qline(-.25,-.25)(-.25,.25)\qline(.25,.25)(-.25,.25)
\qline(-.25,-.25)(.25,-.25)\qline(.25,.25)(.25,-.25)
\psline(-.25,.25)(-.683,.5) \psline(-.25,-.25)(-.683,-.5)
\qline(-.683,.5)(-.683,.75)\qline(-.683,-.5)(-.683,-.75)
\qline(-.683,.5)(-1.116,.25)\qline(-.683,-.5)(-1.116,-.25)
\qline(-1.116,-.25)(-1.116,.25) \psline(-1.116,.25)(-1.549,.5)
\psline(-1.116,-.25)(-1.549,-.5)
\qline(-1.549,.5)(-1.549,.75)\qline(-1.549,-.5)(-1.549,-.75)
\qline(-1.549,.5)(-1.932,.25)\qline(-1.549,-.5)(-1.932,-.25)
\qline(-1.932,-.25)(-1.932,.25) \psline(-1.932,.25)(-2.415,.5)
\psline(-1.932,-.25)(-2.415,-.5)
\psline(.25,.25)(.683,.5)
\psline(.25,-.25)(.683,-.5)
\qline(.683,.5)(.683,.75)\qline(.683,-.5)(.683,-.75)
\qline(.683,.5)(1.116,.25)\qline(.683,-.5)(1.116,-.25)
\qline(1.116,-.25)(1.116,.25) \psline(1.116,.25)(1.333,.375)
\psline(1.116,-.25)(1.333,-.375)
\qline(1.933,.25)(1.717,.375)\qline(1.933,-.25)(1.717,-.375)
\qline(1.933,.25)(1.933,-.25) \psline(1.933,.25)(2.366,.5)
\psline(1.933,-.25)(2.366,-.5)
\qline(2.366,.5)(2.366,.75)\qline(2.366,-.5)(2.366,-.75)
\qline(2.366,.5)(2.799,.25)\qline(2.366,-.5)(2.799,-.25)
\qline(2.799,-.25)(2.799,.25) \psline(2.799,.25)(3.232,.5)
\psline(2.799,-.25)(3.232,-.5)
\end{pspicture}
$$
\caption{Two possible neighborhoods of the moveable
rectangular face.} \label{hexa1}
\end{figure}

\begin{figure}$$
\begin{pspicture}[shift=-2.2](-2.2,-2.2)(1.7,2.2) \rput(1.35,0){\rnode{e1}{$7^+$}}
\rput(0,1.5){\rnode{e3}{$7^+$}} \rput(-1.5,1.5){\rnode{e4}{$6^+$}}
\rput(-2,0){\rnode{e5}{$6^+$}} \rput(-1.5,-1.5){\rnode{e6}{$6^+$}}
\rput(0,-2){\rnode{e7}{$6^+$}} \rput(1.5,-1.5){\rnode{e8}{$6^+$}}
\rput[b](.3,.15){$+$}\rput[b](.3,-.35){$-$}
\rput[b](-.3,.15){$-$}\rput[b](-.3,-.35){$+$}
\rput[b](-1.3,.15){$-$}\rput[b](-1.3,-.35){$+$}
\rput[b](-1,.6){$+$}\rput[b](-1,-.8){$-$}
\rput[b](.3,-1.35){$-$}\rput[b](-.3,-1.35){$+$}
\rput[b](1,-.8){$+$}
\qline(-.5,-.5)(-.5,.5) \qline(-.5,-.5)(.5,-.5)
\qline(.5,.5)(-.5,.5) \qline(.5,.5)(.5,-.5)
\qline(-1.5,-.5)(-1.5,.5) \qline(-.5,-1.5)(.5,-1.5) \qline(1,1)(1,2)
\qline(1,1)(1.5,1) \qline(-1,-1)(-1.5,-.5) \qline(-1,-1)(-.5,-1.5)
\qline(1,-1)(1.5,-1) \qline(1,-1)(.5,-1.5) \qline(-1,1)(-1.5,.5)
\qline(-1,1)(-1,2) \psline(.5,.5)(1,1) \psline(-.5,.5)(-1,1)
\psline(.5,-.5)(1,-1) \psline(-.5,-.5)(-1,-1) \psline(.5,-1.5)(1,-2)
\psline(-.5,-1.5)(-1,-2) \psline(-1.5,.5)(-2,1)
\psline(-1.5,-.5)(-2,-1)
\end{pspicture}
= \sum_i \alpha_i\begin{pspicture}[shift=-2.7](-2.7,-2.7)(2.3,2.2)
\rput(0,0){\rnode{e1}{$7$}}\rput(1.8,.5){\rnode{e2}{$6^+$}}
\rput(.5,2){\rnode{e3}{$6^+$}}\rput(-1.5,1.5){\rnode{e4}{$7^+$}}
\rput(-2.3,0){\rnode{e5}{$6^+$}}\rput(-1.5,-1.5){\rnode{e6}{$7^+$}}
\rput(0,-2.5){\rnode{e7}{$6^+$}}\rput(1.5,-1.5){\rnode{e8}{$7^+$}}
\rput[b](-1.2,.15){$+$}\rput[b](-1.2,-.35){$-$}
\rput[b](-1.8,.15){$-$}\rput[b](-1.8,-.35){$+$}
\rput[b](.3,-1.35){$+$}\rput[b](.3,-1.85){$-$}
\rput[b](-.3,-1.35){$-$}\rput[b](-.3,-1.85){$+$}
\rput[b](.7,.6){$-$}\rput[b](.7,-.4){$+$}
\rput[b](-.35,.6){$+$}
\psline(1,1)(1.5,1.5) \psline(-.5,1)(-1,.5) \psline(-2.5,1)(-2,.5)
\psline(-2.5,-1)(-2,-.5) \psline(-.5,-1)(-1,-.5)
\psline(-.5,-2)(-1,-2.5) \psline(.5,-2)(1,-2.5)
\psline(.5,-1)(1,-.5) \qline(1.5,1.5)(1.5,2)\qline(1.5,1.5)(2,1.5)
\qline(-.5,2.2)(-.5,1) \qline(-.5,1)(1,1) \qline(1,1)(1,-.5)
\qline(1,-.5)(2,-.5) \qline(-1,.5)(-2,.5) \qline(-1,.5)(-1,-.5)
\qline(-2,-.5)(-2,.5) \qline(-2,-.5)(-1,-.5) \qline(.5,-1)(.5,-2)
\qline(.5,-1)(-.5,-1) \qline(-.5,-2)(.5,-2) \qline(-.5,-2)(-.5,-1)
\end{pspicture}  + \sum_j \beta_j ~~~\begin{matrix}\mathrm{webs}~\mathrm{which}\\
\mathrm{have}~ \mathrm{at}~\\
\mathrm{least}~ \mathrm{one}\\\mathrm{less}~
\mathrm{face}\end{matrix}
$$
\caption{Two adjacent hexagons in the neighborhood of the rectangle can not exist
in $D$.
} \label{hexa2}
\end{figure}
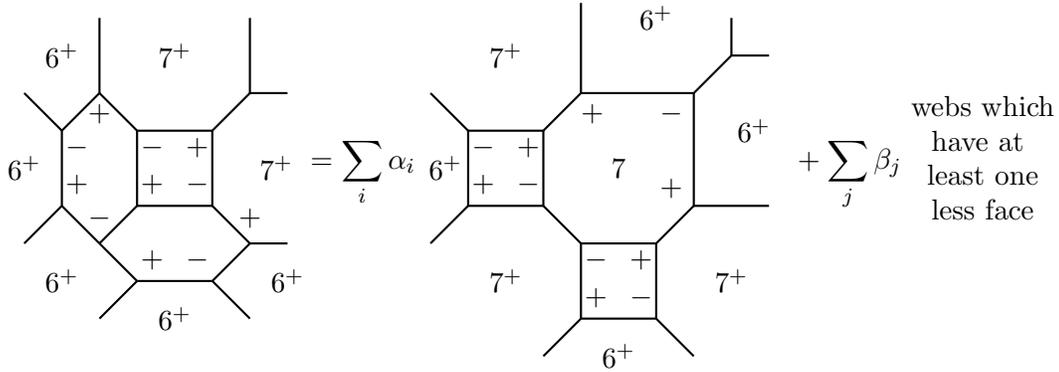
We will divide cases by the number of adjacent
hexagonal faces of the moveable rectangular face.
If there is only one hexagon in the neighborhood of the rectangle, we do swap the rectangle with the hexagon.
After a swapping move, we repeat the process from the new rectangle.
If there is still only one hexagon, then we find the desired result as shown in the left side of
Figure~\ref{hexa1}. If there are more than one hexagon, it will be
dealt with in the next cases.

If there are two adjacent hexagonal faces, there are also two
possibilities. If these two hexagons of valid sign types are adjacent to each other,
we use the relations in Lemma~\ref{rectangular} to change the rectangular
face to a linear combination of webs of less number of faces or webs with
hexagons which can be changed into two rectangles as illustrated in
Figure \ref{hexa2}. Since $D$ is a counterexample,
one of webs in the linear combination has to be a counterexample but
neither cases can be a counterexample because of the
hypothesis of $D$. If these two hexagons are not adjacent each
other, then after a swapping move toward to both directions, we can repeat the process at
the new rectangles. Since there are only finitely many faces, this process either
stops in finite steps which gives a web in the right side of
Figure~\ref{hexa1} or repeats infinitely which gives a circular web obtained from the right one in Figure~\ref{hexa1}
by attaching the rightmost and leftmost hexagons. For three or four adjacent hexagons, we must
have at least two adjacent hexagons which is not possible as described previously.

Furthermore, the traces of swapping moves of different
moveable rectangular faces are disjoint. Otherwise, we can repeat swapping moves to make these two moveable
rectangular faces adjacent. But as we mentioned before, it is not possible to have adjacent rectangles.
\end{proof}

Now, we are set to find a contradictory inequality. Instead of rectangular faces, we are going to look at all polygons
whose sizes are bigger than $6$.
From Lemma~\ref{neighpenta}, three edges from isolated rectangular face and four edges from its adjacent pentagon,
total $4F_5+3(F_4-I)$ edges, can be common edges of faces which have more than seven edges.
From Lemma~\ref{isolatedlem}, for edges in polygons whose size is bigger than $6$ which is
not a common edge with isolated rectangles and their adjacent pentagon, there exists at most one
edge of a moveable rectangular face can travel to become the common edge of the given edge
by a sequence of swapping moves.
Furthermore, since all webs are trivalent and no two rectangles are adjacent, for each face of size $n\ge 9$, there are at most
$ \lfloor\frac{2n}{3}\rfloor$ rectangular faces or pentagonal faces
which are originally adjacent to the given face or adjacent to the given face by a finite sequence of swapping moves.
For a heptagon, it is not adjacent to an isolated rectangle
by Lemma~\ref{neighpenta} and we can see that three rectangles can not be adjacent
to the heptagon because one can easily see that two of these three rectangles can be
adjacent by a swapping move, and having two adjacent rectangles is impossible as mentioned before.
For an octagon, it is not adjacent to an isolated rectangle
by Lemma~\ref{neighpenta} and it can be adjacent to up to $4$ rectangular faces.
These can be summarized as the following inequality,

\begin{align}
4F_5+3(F_4-I)+4I \le 2F_7 + 4F_8 + 6F_9 +6F_{10}+ \ldots + \lfloor
\frac{2n}{3} \rfloor F_n + \ldots.\label{firstinequality}
\end{align}

Since $I=F_4-F_5$, we have $4F_5+3(F_4-I)+4I = 4F_5+3F_4+I=3F_5+4F_4$. Using the fact
$(n-6) \ge \lceil\frac12 \lfloor
\frac{2n}{3} \rfloor\rceil$ for all $n\ge 9$, we
get the following inequalities,

\begin{eqnarray}
4F_4+2F_5&\le& 4F_4+3F_5\le 2F_7+ 4F_8+ 6F_9+ 6F_{10}
+ \ldots + \lfloor \frac{2n}{3} \rfloor F_n + \ldots, \label{secondinequality}\\
 2F_4 +F_5 &\le& 1F_7+2F_8+ 3F_9+4F_{10} + \ldots + (n-6)F_n
+ \ldots, \\
0&\le& -2F_4-F_5+1F_7+2F_8+ 3F_9+4F_{10} + \ldots + (n-6)F_n
+ \ldots . \label{thirdinequality}
\end{eqnarray}

By adding the last inequality~(\ref{thirdinequality}) to the
equality~(\ref{fourthinequality}), we obtain the desired inequality between the number of faces and edges as follows,

\begin{eqnarray}
6F &= 6F_4+6F_5+6F_6+6F_7+6F_8+ \ldots + 6F_n + \ldots,  \label{fourthinequality}\\
6F &\le 4F_4+5F_5+6F_6+7F_7+8F_8+ \ldots + nF_n + \ldots = 2E.  \label{fifthinequality}
\end{eqnarray}

If we substitute the inequality~(\ref{fifthinequality}), $F\le \frac13 E$, and the equation (\ref{graphformula1})
into the equation (\ref{graphformula2}), then we
find a contradiction as

$$2=V-E+F  \le \frac{2}{3}E -E + \frac{1}{3}E =0.$$
Therefore, this completes the proof.
\end{proof}

\section{Proofs of main results}
\label{invariants}

If links are decorated by the fundamental representations
$V_{\lambda_i}$ of the quantum $\mathfrak{sl}(n)$, denoted by $i$,
Murakami, Ohtsuki and Yamada~\cite{MOY:Homfly} found a quantum
invariant for framed links by resolving each crossing in a link
diagram $D$ of $L$ as shown in Figures~\ref{pstvresol1} and
\ref{pstvresol2}.

For negative crossings, we replace $q$ with $q^{-1}$. $P_n(q)$ is
the special case of $q^{-\omega(D)\frac{n}{2}}[D]_n$, when all
components are colored by the fundamental representation
$V_{\lambda_1}$. They showed that $P_n(q)$ is an isotopy invariant
and that $[D]_n$ is a regular isotopy invariant for other colorings
\cite{MOY:Homfly}. However, one can make it a link invariant by
using a suitable writhe but we have to be careful since there are
more than one colors. For a coloring $\mu$ of a diagram $D$ of a
link $L$, we first consider a \emph{colored writhe} $\omega_i(D)$ as
the sum of writhes of components colored by $i$. Then we set
$$K_n(L,\mu) = \prod_i q^{-\omega_i(D)\frac{i(n-i+1)}{2}} [D]_n,$$
where the product runs over all colors $i$.

\begin{figure}
\begin{eqnarray}\left[
\begin{pspicture}[shift=-.6](-.5,-.7)(.5,.7) \rput[bl](.3,.6){$j$}
\rput[br](-.3,.6){$i$} \rput[tr](-.3,-.6){$j$}
\rput[tl](.3,-.6){$i$} \psline(-.25,-.5)(.25,.5)
\psline[arrowscale=1.5]{->}(.15,.3)(.25,.5) \psline(.25,-.5)(.05,-.1)
\psline[arrowscale=1.5]{->}(-.15,.3)(-.25,.5) \psline(-.05,.1)(-.25,.5)
\end{pspicture}\right]_n
=\sum_{k=0}^{i} (-1)^{k+(j+1)i}q^{\frac{(i-k)}{2}} \left[
\begin{pspicture}[shift=-1.2](-1.75,-1.2)(1.75,1.2)
\rput[bl](1.2,.85){$j$} \rput[br](-1.2,.85){$i$}
\rput[tr](-1.2,-.85){$j$} \rput[tl](1.2,-.85){$i$}
\rput[b](0,.5){$j+k-i$} \rput[r](-.9,0){$j+k$}
\rput[t](0,-.5){$k$} \rput[l](.9,0){$i-k$}
\psline(.8,.4)(1.2,.8)\psline[arrowscale=1.5]{->}(.9,.5)(1.1,.7)
\psline(-.8,.4)(-1.2,.8)\psline[arrowscale=1.5]{->}(-.9,.5)(-1.1,.7)
\psline(-1.2,-.8)(-.8,-.4)\psline[arrowscale=1.5]{->}(-1.1,-.7)(-.9,-.5)
\psline(1.2,-.8)(.8,-.4)\psline[arrowscale=1.5]{->}(1.1,-.7)(.9,-.5)
\psline(-.8,.4)(.8,.4)\psline[arrowscale=1.5]{->}(-.1,.4)(.1,.4)
\psline(.8,-.4)(.8,.4)\psline[arrowscale=1.5]{->}(.8,-.1)(.8,.1)
\psline(-.8,-.4)(-.8,.4)\psline[arrowscale=1.5]{->}(-.8,-.1)(-.8,.1)
\psline(.8,-.4)(-.8,-.4)\psline[arrowscale=1.5]{->}(.1,-.4)(-.1,-.4)
\end{pspicture}\right]_n  \label{pstvresol1}
\end{eqnarray}
\begin{eqnarray}\left[
\begin{pspicture}[shift=-.6](-.5,-.7)(.5,.7) \rput[bl](.3,.6){$j$}
\rput[br](-.3,.6){$i$} \rput[tr](-.3,-.6){$j$}
\rput[tl](.3,-.6){$i$} \psline(-.25,-.5)(.25,.5)
\psline[arrowscale=1.5]{->}(.15,.3)(.25,.5) \psline(.25,-.5)(.05,-.1)
\psline[arrowscale=1.5]{->}(-.15,.3)(-.25,.5) \psline(-.05,.1)(-.25,.5)
\end{pspicture}\right]_n
=\sum_{k=0}^{j} (-1)^{k+(i+1)j}q^{\frac{(j-k)}{2}}\left[
\begin{pspicture}[shift=-1.2](-1.75,-1.2)(1.75,1.2)
\rput[bl](1.2,.85){$j$} \rput[br](-1.2,.85){$i$}
\rput[tr](-1.2,-.85){$j$} \rput[tl](1.2,-.85){$i$}
\rput[b](0,.5){$i+k-j$} \rput[r](-.9,0){$j-k$}
\rput[t](0,-.5){$k$} \rput[l](.9,0){$i+k$}
\psline(.8,.4)(1.2,.8)\psline[arrowscale=1.5]{->}(.9,.5)(1.1,.7)
\psline(-.8,.4)(-1.2,.8)\psline[arrowscale=1.5]{->}(-.9,.5)(-1.1,.7)
\psline(-1.2,-.8)(-.8,-.4)\psline[arrowscale=1.5]{->}(-1.1,-.7)(-.9,-.5)
\psline(1.2,-.8)(.8,-.4)\psline[arrowscale=1.5]{->}(1.1,-.7)(.9,-.5)
\psline(-.8,.4)(.8,.4)\psline[arrowscale=1.5]{<-}(-.1,.4)(.1,.4)
\psline(.8,-.4)(.8,.4)\psline[arrowscale=1.5]{->}(.8,-.1)(.8,.1)
\psline(-.8,-.4)(-.8,.4)\psline[arrowscale=1.5]{->}(-.8,-.1)(-.8,.1)
\psline(.8,-.4)(-.8,-.4)\psline[arrowscale=1.5]{<-}(.1,-.4)(-.1,-.4)
\end{pspicture}\right]_n  \label{pstvresol2}
\end{eqnarray}
\caption{Skein expansions of a crossing} \label{expansion}
\end{figure}

\begin{figure}
$$
\left[ \begin{pspicture}[shift=-.6](-.45,-.7)(.55,.7)
\psline(.2,.3)(-.4,-.6)
\psline(-.4,.6)(-.04,.06) \psline(.04,-.06)(.2,-.3)
\psline[arrows=->,arrowscale=1.5](-.1,.15)(-.3,.45)
\pscurve(.2,-.3)(.3,-.4)(.5,0)(.3,.4)(.2,.3)
\end{pspicture}\right]_n = q^{\frac{i(n-i+1)}{2}} \left[
\begin{pspicture}[shift=-.6](-.25,-.7)(0,.7)
\pscurve(-.2,.6)(-.1,.3)(-.05,.1)(-.05,-.1)(-.1,-.3)(-.2,-.6)
\psline[arrows=->,arrowscale=1.5](-.133333,.4)(-.166666,.5)
\end{pspicture} \right]_n ,
\left[ \begin{pspicture}[shift=-.6](-.45,-.7)(.55,.7)
\psline(-.4,.6)(.2,-.3)
\psline[arrows=->,arrowscale=1.5](-.1,.15)(-.3,.45)
\psline(.2,.3)(.04,.06) \psline(-.04,-.06)(-.4,-.6)
\pscurve(.2,-.3)(.3,-.4)(.5,0)(.3,.4)(.2,.3)
\end{pspicture}\right]_n = q^{-\frac{i(n-i+1)}{2}} \left[\begin{pspicture}[shift=-.6](-.25,-.7)(0,.7)
\pscurve(-.2,.6)(-.1,.3)(-.05,.1)(-.05,-.1)(-.1,-.3)(-.2,-.6)
\psline[arrows=->,arrowscale=1.5](-.133333,.4)(-.166666,.5)
\end{pspicture} \right]_n
$$
\caption{Resolving Reidemeister move I.} \label{reider1}
\end{figure}

By Theorem~\ref{onedim}, we know $K_n(L,\mu)$ is well defined. It
was shown that $[D]_n$ is a regular isotopy invariant
\cite{MOY:Homfly}, so is $K_n(L,\mu)$. One can find the equations in
Figure~\ref{reider1} by the second relation in Figure~\ref{firstrel} and a routine
induction. Then $K_n(L,\mu)$ is invariant under Reidemeister move
I. Therefore, $K_n(L,\mu)$ is an isotopy invariant of a link $L$.
A coloring $\mu$ on a $p$-periodic link $L$ is a \emph{$p$-periodic coloring of $L$} if
the periodic homeomorphism $h$ of order $p$ used for the periodicity of $L$
also preserves the coloring $i.e.$, $h(L, \mu)=(L, \mu)$. For such a coloring $\mu$ on a periodic
link $L$, we also denote the factor link $\overline{L}= \pi(L)$ and natural coloring on $\overline{L}$ by
$\overline{\mu}$. Now we discuss the
relation between $K_n(L,\mu)$ and $K_n(\overline{L},\overline{\mu})$ for
a $p-$periodic link $L$ in the following theorems.

\begin{thm}
Let $p$ be a positive integer and $L$ be a $p-$periodic link in
$\mathbb{S}^3$ with its factor link $\overline{L}$. Let $\mu$ be a
$p$-periodic coloring of $L$ and $\overline{\mu}$ be the induced
coloring of $\overline{L}$. Then for $n\ge 0$,

$$K_n(L,\mu) \equiv K_n(\overline{L},\overline{\mu})^p \hskip .7cm modulo \hskip .2cm
\mathcal{I}_n,$$ where $\mathcal{I}_n$ is the ideal of
$\mathbb{Z}[q^{\pm \frac 12}]$ generated by $p$ and
$\left[\begin{matrix}n\\
i\end{matrix}\right]^p-\left[\begin{matrix}n\\
i\end{matrix}\right]$ for $i=1, 2, \ldots, \lfloor \frac{n}{2}
\rfloor$. \label{mainthm}
\end{thm}

\begin{proof}
Let $\mu$ be a $p$-periodic coloring of $L$. Let $R$ be the
fundamental region by the action $h$. Let $\mathcal{C}$ be the set
of all crossings of $L$ and let $\overline{\mathcal{C}}$ be the set
of all crossings in the region $R$. Let $i(c)$ and $j(c)$ be the
weights of two components at the crossing $c$ as in Figure
\ref{expansion}. Let $J(c)$ be the minimum of
$i(c)$ and $j(c)$ for the crossing $c$. Let $D$ be the diagram after
the expansion by the equations in Figures \ref{expansion}. Let $D'$ be the diagram obtained by identical
expanding of the crossings which are the same by the action $h$ and
$D''$ be the diagram obtained from $D'$ by identical applications of
relations at the faces which are in the same orbit by the action
$h$. Let $\overline{D''}= D''/\mathbb{Z}_p$ and $\mathbb{Z}_p$ is
generated by the action $h$. Let $\mathcal{I}_n$ be the ideal of
$\mathbb{Z}[q^{\pm \frac 12}]$ generated by $p$ and
$\left[\begin{matrix}n\\
i\end{matrix}\right]^p-\left[\begin{matrix}n\\
i\end{matrix}\right]$ for $i=1, 2, \ldots, \lfloor \frac{n}{2}
\rfloor$.

\begin{align} \label{eq1}
K_n(L,\mu)&= \prod_{c\in \mathcal{C}}\sum_{k=0}^{J(c)}
(-1)^{k+(i(c)+1)j(c)}q^{\frac{(j(c)-k)}{2}} [D]_n \\
\nonumber &\equiv \prod_{\overline{c}\in
\overline{\mathcal{C}}}\sum_{k=0}^{J(\overline{c})}
((-1)^{k+(i(\overline{c})+1)j(\overline{c})}q^{\frac{(j(\overline{c})-k)}{2}})^p
[D']_n  \hskip 1cm(mod \hskip .3cm p)\\ \nonumber &\equiv
\prod_{\overline{c}\in
\overline{\mathcal{C}}}\sum_{k=0}^{J(\overline{c})}
((-1)^{k+(i(\overline{c})+1)j(\overline{c})}q^{\frac{(j(\overline{c})-k)}{2}})^p
[D'']_n  \hskip 1cm(mod \hskip .3cm p) \\ \nonumber &\equiv
\prod_{\overline{c}\in
\overline{\mathcal{C}}}\sum_{k=0}^{J(\overline{c})}
((-1)^{k+(i(\overline{c})+1)j(\overline{c})}q^{\frac{(j(\overline{c})-k)}{2}})^p
([\overline{D''}]_n )^p \hskip 1cm(mod \hskip .3cm \mathcal{I}_n)
\\ \nonumber
&\equiv  \prod_{\overline{c}\in
\overline{\mathcal{C}}}\sum_{k=0}^{J(\overline{c})}
((-1)^{k+(i(\overline{c})+1)j(\overline{c})}q^{\frac{(j(\overline{c})-k)}{2}} [\overline{D}]_n )^p \hskip 1cm(mod \hskip .3cm p)\\ \nonumber
&\equiv  (\prod_{\overline{c}\in
\overline{\mathcal{C}}}\sum_{k=0}^{J(\overline{c})}
(-1)^{k+(i(\overline{c})+1)j(\overline{c})}q^{\frac{(j(\overline{c})-k)}{2}} [\overline{D}]_n )^p \hskip 1cm(mod \hskip .3cm p)\\ \nonumber
&\equiv  (K_n(\overline{L},\overline{\mu}))^p  \hskip 1cm(mod
\hskip .3cm p)\end{align}

If any expansion of crossings occurs in $R$, it must be used
identically for all other $p-1$ copies of $R$. Otherwise there will
be $p$ identical shapes by the rotation of order $p$, then the term in the expansion is
congruent to zero modulo $p$. This implies the first congruence in
equation (\ref{eq1}). By the same philosophy, if any expansion of
relations occurs in $R$, it must be used identically for all other
$p-1$ copies of $R$. Otherwise it is congruent to zero modulo $p$
and this implies the second congruence. Let us remark that we have
not used relations that might occur in the faces which are not
entirely contained in $R$. If there is an unknot in $\overline{D''}$
which was in $D''$, there are $p$ copies in $D''$. Thus raising
$p$-power gives us the equality. But if there is an unknot in
$\overline{D''}$ which was not in $D''$, it should count once in
$D''$ and $\overline{D''}$, thus we have to use $\mathcal{I}_n$
congruence.
\end{proof}

\begin{rem} In the proof, one can easily  generalize the congruence relation between the quantum invariants $[~~]_n$ of
a web diagram $W$ of symmetry of order $p$ and its quotient $\overline{W}$ modulo $\mathcal{I}_n$ as
$$[W]_n \equiv  ([\overline{W}]_n )^p \hskip 1cm(mod \hskip .3cm \mathcal{I}_n).$$
But in general, this congruence is not true if $\overline{W}$ is a base web diagram~\cite{DL:sl3}.
\end{rem}

Next, we find that Theorem \ref{mainthm} remains true for the colored
$\mathfrak{sl}(n)$ HOMFLY polynomial of a periodic link whose
components are colored by irreducible representations of quantum
$\mathfrak{sl}(n)$.

\begin{thm}
Let $p$ be a positive integer and $L$ be a $p-$periodic link in
$\mathbb{S}^3$ with its factor link $\overline{L}$. Let $\mu$ be a
$p$-periodic coloring of $L$ and $\overline{\mu}$ be the induced
coloring of $\overline{L}$. Then for $n\ge 0$,

$$G_n(L,\mu) \equiv G_n(\overline{L},\overline{\mu})^p \hskip .7cm modulo \hskip .2cm
\mathcal{I}_n,$$ where $\mathcal{I}_n$ is the ideal of
$\mathbb{Z}[q^{\pm \frac 12}]$ generated by $p$ and
$\left[\begin{matrix}n\\
i\end{matrix}\right]^p-\left[\begin{matrix}n\\
i\end{matrix}\right]$ for $i=1, 2, \ldots, \lfloor \frac{n}{2}
\rfloor$. \label{mainthm3}
\end{thm}

\begin{proof}
Since all clasps are idempotents, we put $p$-copies of clasps to
each copy of the fundamental region by the action $h$ of the
periodicity of $L$. Thus without expanding the clasps, we obtain the
theorem by the same idea of the proof of the Theorem~\ref{mainthm}.
\end{proof}

We also give a criterion for periodic links by using the invariant
$K_n(L, \mu)$ and mirror image of knots in the following theorem.

\begin{thm} \label{mainthm2}
\label{mirror2} Let $L$ be a $p$--periodic link for a prime $p$ and
let $L^*$ be the link obtained from the mirror image of a diagram of
$L$. Let $\mu$ be a coloring of $L$ and $\mu^*$ be the coloring of
$L^*$ induced from the coloring $\mu$ of $L$. Then we have
$$K_n(L,\mu) \equiv K_n(L^*,\mu^*)  ~~\mathrm{mod}~~ (p, q^p-1).$$
\end{thm}
\begin{proof}
For a colored link diagram $D$ we denote its mirror image by $D^*$
which is the colored link diagram obtained from $D$ by changing all
of the crossings. We study another necessary condition for a colored
link to be periodic by using the invariant $K_n(L,q)$. Let $L$ be a
colored link diagram with a crossing $x$. Let $L_+,~L_-,~$ and $L_0$
be the link diagrams obtained by resolving the crossing $x$ as shown
in Figure \ref{local} respectively. Then from the Figure
\ref{expansion}, we see that there exist web
diagrams $L_1,\cdots,L_m$ and polynomials $f_1,\cdots,f_m \in
Z[q^{\pm \frac{1}{2}}]$ for some positive integer $m$ such that
$$[L_+]_n - [L_-]_n =
(q^{\frac{1}{2}}-q^{-\frac{1}{2}})(f_1[L_1]_n+\cdots+f_m[L_m]_n).$$

In particular, by applying this relation repeatedly for a colored
periodic link we obtain the following lemma.

\begin{lem} \label{mirror1}
Let $D$ be a colored $p$--periodic link diagram and $\pi$ be the
quotient map of the periodic homeomorphism $h$ of order $p$ and $\overline{D}$ be the factor link of $D$ so
that $\pi^{-1}(\overline{D}) = D$. Let $x$ be a crossing of a
nontrivial diagram $\overline{D}$, and let $\overline{D}_+$ and
$\overline{D}_-$ be the colored link diagram obtained by changing
the crossing $x$ to a positive crossing and negative crossing
respectively. Then we have $$[\pi^{-1}(\overline{D}_+)]_n \equiv
[\pi^{-1}(\overline{D}_-)]_n ~~\mathrm{mod}~~ (p, q^p-1)$$
\end{lem}
\begin{proof}
Let $\overline{D}_1,\cdots,\overline{D}_m$ be web diagrams and let
$f_1,\cdots,f_m \in Z[q^{\pm \frac{1}{2}}]$ be polynomials for some
positive integer $m$ such that
$$[\overline{D}_+]_n - [\overline{D}_-]_n =
(q^{\frac{1}{2}}-q^{-\frac{1}{2}})(f_1[\overline{D}_1]_n+\cdots+f_m[\overline{D}_m]_n).$$

Then by considering the periodic action induced from the map $h$, we
get\\
$[\pi^{-1}(\overline{D}_+)]_n - [\pi^{-1}(\overline{D}_-)]_n$
\begin{eqnarray*}
 &\equiv&
(q^{\frac{1}{2}}-q^{-\frac{1}{2}})^p(f_1^p[\pi^{-1}(\overline{D}_1)]_n
+\cdots+f_m^p[\pi^{-1}(\overline{D}_m)]_n)~~\mathrm{mod}~~ (p)\\
&\equiv& 0 ~~\mathrm{mod}~~ (p, q^p-1).
\end{eqnarray*}
\end{proof}

Now, we are ready to prove theorem. For a tangle $T$, we
denote its closure by $Cl(T)$ if it is well defined. Let $L$ be a
$p$--periodic link and $T$ be a tangle such that $L$ is the closure
$Cl(T^p)$ of the tangle $T^p$ which is $p$ times self-product of
$T$. Let $D$ be a diagram of $L$ and $x$ be a crossing of the
diagram $D$. Let $T_+$ and $T_-$ be the diagram obtained from the
diagram of $T$ by changing the crossing $x$ to a positive crossing
and negative crossing respectively. If the two colorings of the
strands near the crossing $x$ are different then by using Lemma
\ref{mirror1}, we see that

\begin{align*}
&K_n(Cl({(T_+)^p}), \mu) - K_n(Cl({(T_-)^p}), \mu) \\
&= \prod_i
q^{-w_i(D)\frac{i(n-i+1)}{2}}([Cl({(T_+)^p})]_n-[Cl({(T_-)^p})]_n)\\
&\equiv 0 ~~\mathrm{mod}~~ (p, q^p-1).
\end{align*}

If two colorings of the strands near the crossing $x$ are equal, say
$i$, then $w_i(Cl({(T_+)^p})) = w_i(Cl({(T_+)^p})) +2p$. Then we see
that there exists a polynomial $g\in \mathbb{Z}[q^{\pm \frac{1}{2}}]$
such that

\begin{align*}
&K_n(Cl({(T_+)^p}), \mu) - K_n(Cl({(T_-)^p}), \mu) \\
&=
g([Cl({(T_+)^p})]_n-q^{pi(n-i+1)}[Cl({(T_-)^p})]_n)\\
&\equiv 0 ~~\mathrm{mod}~~ (p, q^p-1).
\end{align*}

The last congruence relation in the above formulae can be obtained
by using Lemma \ref{mirror1}. Thus we see that

$$K_n(Cl({(T_+)^p}), \mu) \equiv K_n(Cl({(T_-)^p}), \mu)\quad\mathrm{mod}~~ (p, q^p-1).$$

Since $T$ can be obtained from its mirror image $T^*$ by a finite
sequence of crossing changes, we get
$$K_n(L,\mu) \equiv K_n(L^*,\mu^*) \quad\mathrm{mod}~~ (p, q^p-1).$$

It completes the proof the theorem.
\end{proof}

\section{Discussion}\label{discussion}

To compare our result with previous results,
one has to compare the size of ideal used in Theorem~\ref{mainthm} with ideals in~\cite{chbili, CL:period, PS:superiod, mu:jones}.
As we have mentioned, we have corrected the false conjecture given by Chbili~\cite{chbili} and Przytycki and Sikora
~\cite{PS:superiod}. To compare with the ideal in~\cite{mu:jones}, we observe that
the ideal $\mathcal{I}_n$ of
$\mathbb{Z}[q^{\pm \frac 12}]$ generated by $p$ and
$\left[\begin{matrix}n\\
i\end{matrix}\right]^p-\left[\begin{matrix}n\\
i\end{matrix}\right]$ for $i=1, 2, \ldots, \lfloor \frac{n}{2}
\rfloor$ is a {\bf subset} of the ideal generated by $p$ and $[2]^p-[2]$
by the strong integrality of the quantum link invariant~\cite{Le:integral}.
Furthermore, if $n$ is odd, the ideal $\mathcal{I}_n$ is a subset of
the ideal generated by $p$ and $[3]^p-[3]$. To compare the ideal generated by $p$ and $[2]^p-[2]$ with the ideal of Murasugi's~\cite{mu:jones} generated by $p$ and
$$\xi_p(t)=\sum_{j=0}^{p-1} (-t)^j-t^{\frac{p-1}{2}},$$
we observe
$$(t+1)\xi_p(t) \equiv q^{-\frac{p}2}([2]^p-[2])|_{\sqrt{q}=-\frac{1}{\sqrt{t}}} ~~~~(\mathrm{mod}~p).$$
To compare with the ideal in~\cite{CL:period}, we use only fundamental representations of the quantum Lie algebras $\mathfrak{sl}(n)$ which are finite but Chen and Le used all representations of the quantum $\mathfrak{sl}(n)$ which are obviously infinite.
Thus, our criterion is sharper than other previous results.
However, we were not able to find a new periodicity of knots using our criteria.
\vskip 1cm
 \noindent{\bf Acknowledgements}

The authors would like to thank Greg Kuperberg for introducing the
subject and helpful discussion, Younghae Do, Mikhail Khovanov, Hitoshi
Murakami for their attentions to this work. Also, the referee has
been very critical during refereing and pointing out the comparison in Section~\ref{discussion}. The \TeX\, macro
package PSTricks~\cite{PSTricks} were essential for typesetting the
equations and figures.




\begin{thebibliography}{00}

\bibitem{Adams} C. Adams, The Knot Book: An Elementary Introduction to the
Mathematical Theory of Knots. New York: W. H. Freeman, 1994.

\bibitem{BHMV} C. Blanchet, N. Habegger, G. Masbaum and P. Vogel, \textit{Topological
Quantum Field Theories derived from the Kauffman bracket}, Topology
34 (1995), 883--927.

\bibitem{CK:homology} S. Cautis and J. Kamnitzer, \textit{Knot homology via derived categories of coherent sheaves I,
SL(2) case}, Duke Math. J. Volume 142, Number 3 (2008), 511--588.

\bibitem{chbili} N. Chbili, \textit{The quantum SU(3) invariant of links and Murasugi's congruence},
Topology appl., 122, (2002), 479--485.

\bibitem{chbili:qm} N. Chbili, \textit{Quantum invariants and finite group actions on three-manifolds}, Topology appl., 136, (2004), 219--231.

\bibitem{CL:period} Q. Chen and T. Le, \textit{Quantum invariants and periodic links and periodic manifolds}, Fund.
Math. 184 (2004), 55--71.

\bibitem{FK:canonical} I. Frenkel and M. Khovanov, \textit{Canonical bases in tensor products and
graphical calculus for $U_{q}(\mathfrak{sl}_2)$}, Duke Math. J.,
87(3), (1997) 409--480.

\bibitem{FultonHarris:gtm} W. Fulton and J. Harris, \textit{Representation theory},
Graduate Texts in Mathematics, 129, Springer-Verlag, New
York-Heidelberg-Berlin, 1991.

\bibitem{Jones:subfactor} V. F. R. Jones, \textit{Index of subfactors}, Invent.
Math.,  72 (1983), 1--25.

\bibitem{Jones:braid} V. F. R. Jones, \textit{Hecke algebra representations of braid groups and link polynomials},
Ann. of Math., 126 (1987), 335--388.

\bibitem{JP5} M.--J. Jeong and C.--Y. Park, \textit{Lens knots, periodic knots and Vassiliev invariants}, J. of Knot Theory and Its
Ramifications, Vol. 13 (2004), 1041--1056.

\bibitem{KRT:knot} C. Kassel, M. Rosso and V. Turaev, \textit{Quantum groups and knot invariants},
Panoramas et Syntheses, 5, Societe Mathematique de France, 1997.

\bibitem{Khovanov:sl3} M. Khovanov, \textit{$sl(3)$ link homology}, Algebr. Geom. Topol.,
4 (2004), 1045--1081.

\bibitem{Khovanov:colored} M. Khovanov, \textit{Categorifications of the colored Jones polynomial},
J. Knot Theory Ramifications, 14(1) (2005), 111--130.

\bibitem{khvanov:comm} M. Khovanov, private communication.

\bibitem{KR:factor} M. Khovanov and L. Rozansky, \textit{Matrix factorizations and link homology}, Fund. Math., 199 (2008), 1--91.

\bibitem{Dongseok:thesis} D. Kim, \textit{Graphical Calculus on Representations of Quantum Lie
Algebras}, Thesis, UCDavis, 2003, arXiv:math.QA/0310143.

\bibitem{DL:sl3} D. Kim and J. Lee, \textit{The quantum sl(3) invariants of cubic bipartite planar graphs},
J. Knot Theory Ramifications, 17(3) (2008), 361--375.

\bibitem{kirbyMelvin:witten} R. Kirby and P. Melvin,  \textit{The 3-manifold
invariants of Witten and Reshetikhin-Turaev for $\mathfrak{sl}(2)$},
Invent. Math., 105 (1991), 473--545.

\bibitem{Kuperberg:spiders} G. Kuperberg, \textit{Spiders for rank 2 {Lie} algebras}, Comm. Math. Phys.,
180(1), (1996) 109--151.

\bibitem{Le:integral} T. Le, \textit{Integrality and symmetry of quantum link invariants},
Duke Math. J., 102 (2000), 273--306.

\bibitem{Lickorish:su} W. Lickorish, \textit{Distinct 3-manifolds with all $SU(2)_q$ invariants the same},
Proc. Amer. Math. Soc., 117 (1993), 285--292.

\bibitem{Morrison} Scott Morrison,
\textit{A Diagrammatic Category for the Representation Theory of $U_q(sl_n)$}, UC Berkeley Ph.D. thesis, arXiv:0704.1503.

\bibitem{mu:alexander} K. Murasugi, \textit{On periodic knots}, Comment. Math. Helv.,
46 (1971), 162--174.

\bibitem{mu:jones} K. Murasugi, \textit{The Jones polynomials of periodic links},
Pacific J. Math., 131 (1988), 319--329.

\bibitem{Murakami:coloredjones} H. Murakami, \textit{Asymptotic Behaviors of the colored
Jones polynomials of a torus knot}, Internat. J. Math.,  15(6) (2004), 547--555.

\bibitem{MOY:Homfly} H. Murakami and T. Ohtsuki and S. Yamada, \textit{HOMFLY polynomial via an invariant of colored plane graphs},
L'Enseignement Mathematique, t., 44 (1998), 325--360.

\bibitem{OY:quantum} T. Ohtsuki and S. Yamada, \textit{Quantum $su(3)$ invariants via linear skein theory},
J. Knot Theory Ramifications,  6(3) (1997), 373--404.

\bibitem{Przytycki:criterion} J. H. Przytycki, \textit{On Murasugi's and Traczyk's criteria for
periodic links}, Math. Ann. 283 (1989), 465--478.

\bibitem{PS:superiod} J. Przytycki and A. Sikora, \textit{$SU_n$-Quantum Invariants for Periodic Links, Diagrammatic morphisms and applications}, Contemp. Math., 318 (2003), 199--205.

\bibitem{PS:skein} J. Przytycki and A. Sikora, \textit{On skein algebras and Sl2(C)-
character varieties}, Topology 39 (2000), 115--148.

\bibitem{RT:1} N. Yu. Reshetikhin and V. G. Turaev,
\textit{Ribbob graphs and their invariants derived from quantum
groups}, Comm. Math. Phys., 127 (1990), 1--26.

\bibitem{RT:2} N. Yu. Reshetikhin and V. G. Turaev,
\textit{Invariants of $3$-manifolds via link polynomials and quantum
groups}, Invent. Math., 103 (1991), 547--597.

\bibitem{Turaev:skein} Turaev V. G.,textit{The Conway and Kauffman modules of a solid torusa},
(translation) J. Soviet Math. 52(1) (1990), 2799--2805.

\bibitem{SW:graph}
A. Sikora, B. Westbury, \textit{Confluence theory for graphs}, Algebraic \& Geometric Topology, 7 (2007), 439--478.

\bibitem{traczyk:period3} P. Traczyk, \textit{A criterion for knots of period $3$}, Topology and its Appl. 36 (1990), 275--281.

\bibitem{PSTricks}
T.~{Van Zandt}. PSTricks: {PostScript} macros for generic {\TeX}.
Available at {\tt ftp://ftp.princeton.edu/ pub/tvz/}.

\bibitem{Vybornov:yang} M. Vybornov, \textit{Solutions of the Yang-Baxter equation
and quantum $\mathfrak{sl} (2)$}, J. Knot Theory Ramifications,
8(7) (1999), 953--961.

\bibitem{Wenzl:Proj} H. Wenzl, \textit{On sequences of projections},
C. R. Math. Rep. Acad. Sci. R. Can., IX (1987), 5--9.

\bibitem{math.QA/0601209}
B. Westbury, \textit{Invariant tensors for the spin representation
of $\mathfrak{so}(7)$}, Math. Proc. Cam. Phil. Soc., 144(1) (2008), 217--240.

\bibitem{Witten:pathint} E. Witten, \textit{Quantum field theory and the Jones polynomial},
Commun. Math. Phys., 121 (1989), 300--379.

\bibitem{yokota:skein} Y. Yokota, \textit{The skein polynomial of periodic knots}, Math. Ann.
291(2) (1991), 281--291.

\bibitem{yokota:Jones} Y. Yokota, \textit{The Jones polynomial of periodic knots}, Proc. Amer. Math. Soc.,
113(3) (1991), 889--894.

\bibitem{yokota:kauffman} Y. Yokota, \textit{The Kauffman polynomial of periodic knots}, Topology
32(2) (1993), 309--324.

\bibitem{yokota:skeinforn} Y. Yokota, \textit{Skein and quantum
$SU(N)$ invariants of 3-manifolds}, Math. Ann., 307 (1997),
109--138.

\end{thebibliography}
\end{document}